\documentclass[column]{autart}    
\usepackage{amsmath}
\usepackage{graphicx}
\usepackage{subfigure}
\usepackage{amsmath}
\usepackage{algorithm}
\usepackage{algpseudocode}
\usepackage{cite}
\usepackage{multicol}
\usepackage{multirow}
\usepackage{mathtools}
\usepackage{amssymb}
\usepackage{color}
\usepackage{indentfirst}
\usepackage[hyphens]{url}
\usepackage[hidelinks]{hyperref}
\usepackage{nomencl}
\usepackage{booktabs}
\usepackage{bm}
\usepackage[utf8]{inputenc}
\usepackage[numbers,sort&compress]{natbib}
\usepackage{makecell}
\usepackage{scalerel}
\usepackage{caption}
\usepackage{siunitx}
\newcommand{\smallsum}{{\textstyle \sum}}

\usepackage[hidelinks]{hyperref}

\newtheorem{theorem}{Theorem}

\newtheorem{proposition}{Proposition}
\newtheorem{corollary}{Corollary}
\newtheorem{definition}{Definition}
\newtheorem{remark}{Remark}

\newcommand{\Lag}{\mathbb{L}}
\newcommand{\x}{\mathbf{x}}
\newcommand{\starx}{\mathbf{x}^{*}}
\newcommand{\hx}{\bar{\mathbf{x}}}

\newcommand{\w}{\mathbf{w}}
\newcommand{\bb}{\mathbf{b}}

\newcommand{\y}{\mathbf{y}}
\newcommand{\X}{\mathbf{X}}
\newcommand{\A}{\mathbf{A}}
\newcommand{\B}{\mathbf{B}}
\newcommand{\Q}{\mathbf{Q}}

\newcommand{\R}{\mathbf{R}}
\newcommand{\T}{\mathbf{T}}
\newcommand{\K}{\mathbf{K}}

\newcommand{\m}{\mathbf{m}}

\newcommand{\boldlambda}{\boldsymbol{\lambda}}

\DeclareMathOperator{\Lips}{\text{Lipschitz}}

\begin{document}

\begin{frontmatter}

\title{ Proximal ADMM for Nonconvex   and Nonsmooth Optimization} 

\thanks[footnoteinfo]{
	This  work  is  supported  by  the  Republic  of  Singapore’s  National  Research  Foundation  through  a  grant  to  the  Berkeley  Education  Alliance  for  Research  in  Singapore
	(BEARS)  for  the  Singapore-Berkeley  Building  Efficiency  and  Sustainability  in  the
	Tropics  (SinBerBEST)  Program. This work is  also supported in part by National Natural Science Foundation of China (62192752, 62192750, 62125304, 62073182), 111 International Collaboration Project (BP2018006), and Tsinghua University Initiative Scientific Research Program. \\
\texttt{Yu Yang is the corresponding author. } }

\author[XJTU]{Yu~Yang}\ead{yangyu21@xjtu.edu.cn},   
\author[TS]{Qing-Shan~Jia}\ead{jiaqs@tsinghua.edu.cn},   
\author[XJTU]{Zhanbo~Xu}\ead{zhanbo.xu@xjtu.edu.cn},    
\author[XJTU,TS]{Xiaohong~Guan}\ead{xhguan@xjtu.edu.cn},            
\author[Berkeley]{ Costas J.~Spanos}\ead{spanos@berkeley.edu}

\address[XJTU]{School of Automation Science and Engineering, Xi’an Jiaotong University, Shaanxi, China.}  
\address[TS]{CFINS, Department of Automation, BNRist, Tsinghua University, Beijing, China.}                                                  
\address[Berkeley]{Electrical Engineering and Computer Sciences, University of California, Berkeley.}


\begin{abstract}                          
By enabling the nodes or agents to solve small-sized subproblems to achieve coordination, distributed algorithms are favored  by many networked systems for efficient and scalable computation.  
While for convex problems, substantial distributed algorithms are available,  the results for the more broad nonconvex counterparts are extremely lacking.  
This paper develops a distributed algorithm for  a class of nonconvex and nonsmooth problems 
featured by i)  a nonconvex  objective formed by  both separate and composite components regarding the decision variables of  interconnected agents,  ii) local bounded convex constraints,  and iii) coupled  linear constraints. This problem is directly originated from  smart buildings and is also broad in other domains. To provide a distributed algorithm with convergence guarantee, we revise the powerful   alternating direction method of multiplier (ADMM) method and proposed a proximal ADMM.  Specifically, noting that the main difficulty to establish the convergence  for the nonconvex and nonsmooth optimization with  ADMM is to assume the boundness  of dual updates, we propose to update  the dual variables in a discounted manner.  This leads to the establishment of a so-called 
sufficiently decreasing and lower bounded Lyapunov function, which is critical to establish the convergence. 
We prove that the method converges to some approximate stationary points.
We besides showcase the efficacy and performance of the method by a numerical example and the concrete application to multi-zone heating, ventilation, and air-conditioning (HVAC) control  in smart buildings. 

\end{abstract}

\end{frontmatter}

\section{Introduction}
By enabling the nodes or agents to solve small-sized subproblems to achieve coordination, distributed algorithms are favored by many networked systems to achieve efficient and scalable computation. 
While distributed algorithms for convex optimization  have been studied extensively  \cite{shi2014linear, deng2017parallel,  falsone2020tracking}, 
the results for the more broad  nonconvex counterparts are extremely lacking. 
The direct extension of distributed algorithms for convex problems to nonconvex counterparts  is in general not applicable either due to the failure of convergence 
or the lack of theoretical convergence guarantee ( see \cite{houska2016augmented, wang2019global} for some divergent examples).
This paper focuses on developing  a distributed algorithm for a class of nonconvex and nonsmooth problems in the canonical form of
\begin{subequations}  
	\setlength{\abovedisplayskip}{-2pt}
	\setlength{\belowdisplayskip}{1pt}
	\begin{align}
		\label{pp:problem template} \tag{$\mathbf{P}$}  \min_{\x =(\x_i)_{i=1}^N } &F(\x)= g(\x) +  \sum_{i=1}^N f_i(\x_i) \\
		\label{eq:1a} \tag{1a} {\rm s. t.}~ & \sum_{i=1}^N \A_i \x_i = \bb. \\
		\label{eq:1b} \tag{1b}& \x_i \in \X_i,  ~i = 1, 2, \cdots, N. 
		\end{align}
\end{subequations}
 where $i = 1, 2, \cdots, N$ denotes the computing nodes or agents,  $\x_i \in \R^{n_i}$ is the local decision variables of agent $i$ and $\x = (\x_{i})_{i = 1}^N \in \R^n$ with $n={\textstyle \sum}_{i=1}^N n_i$ stacks the decision variables of all agents.  We have  $f_i: \R^{n_i} \rightarrow \R$ and $g: \R^{n} \rightarrow \R$ denote the separate and composite objective components,  which are continuously differentiable but possibly nonconvex. We have $\X_i$ represent the local bounded and convex constraints of agent $i$. 
As expressed by the formulation,  the agents are expected to  optimize their  local decision variables in a cooperative manner  so as to 
achieve the optimal system performance measured by   $F(\x) = g(\x) + {\textstyle \sum}_{i = 1}^N f_i(\x_i)$ considering both their local  constraints $\X_i$ and the global coupled linear constraints \eqref{eq:1a}  encoded by $\A_i \in \R^{m \times n_i}$ and $b \in \R^{m}$.
By defining $\A = (\A_1, \A_2, \cdots,  \A_N) \in \R^{m \times n}$ and $f(\x) = \sum_{i=1}^N f_i(\x_i)$, the coupled  constraints and objective  can be expressed by  $\A \x = \bb$ and $F(\x) = f(\x) + g(\x)$.  Note that the presence of local constraints 
$\X_i$  and  nonconvex objectives $f_i$ and $g$  makes the problem nonconvex and nonsmooth, which represents the major challenge to develop distributed algorithm with convergence guarantee. 


Problem \eqref{pp:problem template} is directly originated from smart buildings where smart devices are empowered to make local decisions while accounting for the interactions or the shared resource limits with  the other devices in the proximity  (see, for examples \cite{yang2020hvac, yang2021distributed}).    
Many other applications also fit into this formulation, including but not limited to
smart sensing \cite{ansere2020optimal}, energy storage sharing \cite{yang2021optimal, yang2020selling}, electric vehicle charging management \cite{zhang2016scalable, yang2018decentralized, yang2017distributed, yang2017stochastic, yang2016joint, long2021efficient}, peer-to-peer energy trading \cite{yang2022optimal, chen2022towards}, power system control \cite{arpanahi2020comprehensive},  wireless communication control \cite{hashempour2021distributed}. 
When the number of nodes is large, centralized methods usually suffer bottlenecks from the heavy computation,  data storing and  communication (see \cite{yang2020hvac, arpanahi2020comprehensive, hashempour2021distributed} and the references therein).  Also, centralized methods may disrupt privacy as the complete information of  all agents (e.g., the private local objectives) are required by a central computing agent.   
As a result, distributed algorithms are usually preferred for  privacy, computing efficiency, small data storage, and scaling properties.

When problem \eqref{pp:problem template} is convex, plentiful distributed solution methods are available. The methods can be distinguished by the presence of the composite objective component $g$ and the number of decision blocks $N$. When $g$ is null,  we have the classic dual decomposition methods  \cite{necoara2015linear, falsone2017dual}, the well-known alternating direction method of multiplier (ADMM)  for two decision blocks ($N=2$) \cite{boyd2011distributed} and the variations for  multi-block settings  ($N \geq 3$)\cite{lin2015sublinear, cai2017convergence, bai2018generalized}.  
While the classic ADMM and its variations propose to update the decision components in a sequential manner (usually called \emph{Gauss-Seidel} decomposition),  the works  \cite{deng2017parallel} and \cite{chatzipanagiotis2017convergence}  have made some effort in developing parallel ADMM and its variations  (usually called  \emph{Jacobian ADMM} or \emph{parallel ADMM}).  The above methods are generally limited to  separable objective functions (i.e., only $f_i$ exist and $g=0$). 
For the general case with composite objective component $g$,  linearized ADMM \cite{aybat2017distributed} and inexact linearized ADMM \cite{bai2022inexact} are also studied.

The above results are all for convex problems.
Nevertheless, massive applications arising from the engineering systems and machine learning  domains require to handle the type of problem \eqref{pp:problem template} with possibly nonconvex objectives $f_i$ and $g$.  The non-convexity may originate from the complex system performance metrics or the penalties  imposed on the operation constraints. 
When the objectives $f_i$  and $g$ lack convexity (i.e.,  the monotonically non-decreasing property of gradients or subgradients is lost),   developing distributed methods with theoretical convergence guarantee  becomes a  much more challenging problem.  
Though some fresh distributed methods for constrained nonconvex problems  have been developed, they can not be applied to problem \eqref{pp:problem template} due to the  nonsmoothness caused by the local constraints $\X_i$. This can be perceived from the following literature.



\begin{table*}[h] 
	\setlength\tabcolsep{0.6pt}
	\centering
	\caption{Distributed constrained nonconvex optimization}
	\label{tab:literature}
	\begin{tabular}{clccccc}     
		\toprule
		\textbf{\# Type} & \textbf{~Problem structures} &  \textbf{Main assumptions}  &\textbf{Methods}   & \textbf{Scheme} & \textbf{Convergence} & \textbf{Papers}\\
		\hline 
1	&	$\begin{aligned}
			&\min_{(\x_i)_{i=1}^N} \sum_{i=1}^N f_i(\x_i) \\
			& {\rm s. t.}~  \sum_{i=1}^N \A_i \x_i = \bb. \\
			& \x_i \in \X_i,  ~i = 1, 2, \cdots, N.
		\end{aligned}$ 
		&  \makecell{$f_i$ continuously\\ differentiable.  \\ Strong second-order \\ optimality condition.}
		&   ADAL  &Jacobian  & \makecell{Local convergence. \\ Local optima.} & \cite{chatzipanagiotis2017convergence}\\
		\hline
	2 &	$\begin{aligned}
			& \min_{\x =(\x_i)_{i=0}^p,  y} \!\!\!\!g(\x) \!+\! \sum_{i=0}^p f_i(\x_i) \!+\! h(\y)  \\
			& {\rm s. t.}~  \sum_{i=0}^p \A_i \x_i + \B \y = 0.  \\
		\end{aligned}$ 
		&  \makecell{$g$ and $h$ Lipschitz \\ continuous gradient.  \\  $f_i$ weakly convex.   \\
			$\text{Im}(\A) \subseteq \text{Im}(\B)$.  }	
		&   ADMM & Gauss-Seidel  &  \makecell{Global convergence.  \\ Stationary points.} & \makecell{ \cite{wang2019global, yang2017alternating} \\ \cite{guo2017convergence, li2015global} }\\
		\hline
	3 &			$\begin{aligned}
					& \min_{\x=(\x_i)_{i=1}^N,  y} \!\!\!\!\!\! g(\x, \y) \!\!+\!\! \sum_{i=1}^N \!f_i(\x_i) \!+\! h(\y)  \\
					& {\rm s. t.}~  \sum_{i=1}^N \A_i \x_i + \B \y = 0.  \\
				\end{aligned}$ 
				&  \makecell{$g$ and $h$ Lipschitz \\ continuous gradient.  \\
					$\text{Im}(\A) \subseteq \text{Im}(\B)$.  }	
				&  \makecell{ Linearized \\ ADMM} &\makecell{ Gauss-Seidel~~~}  &  \makecell{Global convergence.  \\ Stationary points.} &\cite{liu2019linearized}\\
				\hline
		4 & 		$\begin{aligned}
			&\min_{(\x_k)_{k=0}^K} \sum_{k=1}^{k} g_{k}(\x_{k}) + h(\x_{0})\\
			 &{\rm s. t.}~  \x_{k} = \x_{0}. \\
			& \x_0 \in \X.
		\end{aligned}$ 
		&  \makecell{$g$ Lipschitz \\ continuous gradient. \\  $h$ convex. }	
		&  \makecell{Flexible\\ ADMM} & Gauss-Seidel  &  \makecell{Global convergence.  \\ Stationary points.} & \cite{hong2016convergence}\\
				\hline
	5 &	$\begin{aligned}
			& \min_{ (\x_k)_{k=0}^K} \sum_{k=1}^N g_x(\x_k) + \ell(\x_0) \\
			&{\rm s. t.} ~ \sum_{k=1}^K \A_k \x_k = \x_0. \\
			& \x_k \in \X_k,  ~k = 1, \cdots, N.
		\end{aligned}$ 
		&  \makecell{$\ell$ $\Lips$ \\
			continuous gradient. \\
			$g$ nonconvex but smooth\\
			or convex but non-smooth.}
		&   \makecell{Flexible\\ ADMM} &Gauss-Seidel   & \makecell{Global convergence.  \\
			Stationary points.}  & \cite{hong2016convergence}\\		
		\hline
	6 &	$\begin{aligned}
			& \min_{ (\x_i)_{i=1}^N,  \hx}\sum_{i=1}^N f_i(\x_i) \\
			& {\rm s. t.}~ \sum_{i=1}^N \A_i \x_i + \B \hx = 0. \\
			& \x_i \in \X_i,  h_i(\x_i) = 0, \\
			& \quad \quad   i = 1, \cdots, N. \\
			&\hx \in \bar{\mathbf{X}}.
		\end{aligned}$ 
		&  \makecell{$f_i$ continuously \\
			 differentiable.  \\ $h_i$ non-linear \\
			(possibly nonconvex). \\ $\B$ full column rank. \\
			$\X_i$ possibly nonconvex. }	
		&  \makecell{ ALM \!+\\\! ADMM} & Gauss-Seidel  & \makecell{Global convergence.  \\ Stationary points.} & \cite{sun2019two, sun2021two, yang2020proximal}\\
		\hline	
	7 & 	$\begin{aligned}
			& \min_{ (\x_i)_{i=1}^N} g(\x) + \sum_{i=1}^N f_i(\x_i) \\
			& {\rm s. t.} ~ \sum_{i=1}^N \A_i \x_i = \bb. \\
			& \x_i \in \X_i,  ~i = 1, 2, \cdots, N.
		\end{aligned}$ 
		&  \makecell{$f_i$ and $g$ Lipschitz \\ continuous gradient.}
		&   \makecell{Proximal \\ ADMM} &Jacobian  & \makecell{Global convergence. \\
			Approximate \\ stationary points.}  & This paper\\
		\bottomrule 
	\end{tabular}
\begin{tabular}{l}
Note: the set $\X_i$ and $\bar{\mathbf{X}}$ are bounded convex sets. ~~~~~~~~~~~~~~~~~~~~~~~~~~~~~~~~~~~~~~~~~~~~~~~~~~~~~~~~~~~~~~~~~~~~~~~~~~~~~~~~~~~~~~~~~~~~~~~~~~~~~~~~~~
\end{tabular}
\end{table*}

A comprehensive survey of ADMM for constrainted optimization is available \cite{yang2022survey}.  
The existing works  for constrained nonconvex optimization  can be distinguished by \texttt{problem structures}, \texttt{main assumptions}, \texttt{decomposition scheme} (i.e., \emph{Jacobian} or \emph{Gaussian-Seidel})  and \texttt{convergence guarantee} as reported in Table \ref{tab:literature}.  Overall,  they can be uniformly expressed by the template of problem \eqref{pp:problem template} but are slightly different in the settings and assumptions. 


The first category (Type 1) is concerned with problem \eqref{pp:problem template}   without any composite objective component $g$  \cite{chatzipanagiotis2017convergence}.  
An accelerated distributed augmented Lagrangian (ADAL) method  was  proposed to handle the possibly nonconvex but continuously differentiable  objectives $f_i$.  
This method follows the classic ADMM framework  but introduces an  interpolation procedure  regarding the primal updates  at each iteration,  which reads as $\A_i \x_i^{k+1} = \A_i \x_i^k + \mathbf{T}\left( \A_i  \hat{\x}^k_i - \A_i \x^k_i\right) $ ($k$ the iteration  and $\T$ is a weighted matrix).  
To our understanding, this can be interpreted as a means to  slow down  the primal update for enhancing the convergence in nonconvex settings. 
By assuming  the existence of stationary points that satisfy the strong second-order  optimality condition, this paper established the local convergence  of the method.  
The notion of local convergence is that the convergence towards some local optima can be assured if  starting  with a point  sufficiently close to that local optima.

The subsequent four categories (Type 2, 3, 4, 5) differ from the first  one mainly in  the presence of a last block encoded by $\B$.  Note that  \cite{hong2016convergence}  can be viewed as a special case with $\B = \mathbf{I}$,  where $\mathbf{I}$ are identity matrices of suitable sizes.
The last block is  exceptional  due to the unconstrained and   $\Lips$ differentiable property,  which are critical to bound the dual updates  for establishing convergence (see the references therein).    That's why the last decision block is usually distinguished by some special notations (i.e., $\y$, $\x_0$).  
While the first category employs \emph{Jacobian} decomposition for primal update,  these four categories fall into \emph{Gauss-Seidel} decomposition (i.e., alternating minimization).  Specially,  the works \cite{wang2019global}  and \cite{liu2019linearized} have made some effort in handling possible composite objective components $g$ but via different ways. Specifically, \cite{wang2019global} employed block coordinate and \cite{liu2019linearized}  used  linearization technique. 
Particularly,  \cite{hong2016convergence, wang2019global} build a general framework  to establish the convergence for  \emph{Gauss-Seidel} ADMM towards local optima or stationary points in nonconvex settings, which comprises two key steps: 1)  identifying  a so-called sufficiently  decreasing Lyapunov function, and 2) establishing  the lower boundness property of the Lyapunov function.  The  sufficiently decreasing  and lower boundness property of a proper Lyapunov function state that  \cite{wang2019global}
\begin{equation}\label{eq:sufficiently-decreasing}
	\begin{split}
	\setlength{\abovedisplayskip}{-2pt}
	\setlength{\belowdisplayskip}{1pt}
&T(\x^{k+1}, \boldlambda^{k + 1}) -  T(\x^k, \boldlambda^k) \\
 & \quad \quad \quad \leq  - a_{\x} \Vert \x^{k +1} - \x^{k} \Vert^2 - a_{\boldlambda} \Vert \boldlambda^{k +1} - \boldlambda^k \Vert^2.  \\
&T(\x^k, \boldlambda^k)  > -\infty. 
\end{split}
\end{equation}
where  $T(\cdot, \cdot)$ is a general Lyapunov function, $\x$ and $\boldlambda$ are primal and dual variables,  $a_{\x}$ and $a_{\boldlambda}$ are positive coefficients.
The augmented Lagrangian (AL) function has been often used as the Lyapunov function in nonconvex settings  (see \cite{hong2016convergence, wang2019global} and the references therein).  However, they  depend on  the   following two necessary conditions on  the last decision block encoded by $\B$ to bound the dual updates $\Vert \boldlambda^{k+1} - \boldlambda^k \Vert^2$ by the primal updates $\Vert \x^{k+1} - \x^k \Vert^2$ \cite{hong2016convergence, wang2019global}. 
\begin{itemize}
	\item [a)] $\B$ has full column rank and $\text{Im}(\A) \subseteq \text{Im}(\B)$ ($\text{Im}(\cdot)$ represents the image of a matrix).
	\item [b)]  The last decision block is unconstrained and with $\Lips$  differentiable objective.  
\end{itemize}
Noted that  the forth and fifth category (Type 4, 5) originated from \cite{hong2016convergence}  are  a special case with $\B =\mathbf{I}$ and  thus satisfy  the necessary \textbf{condition} a).

Following the line of works, the sixth category (Type 6)  studied  the extension of ADMM to non-linearly constrained nonconvex problems \cite{sun2019two, sun2021two, yang2020proximal}.  Since it is difficult (if not impossible) to directly handle the non-linear couplings  by the AL framework, \cite{sun2019two} proposed to  first convert the non-linearly constrained problems to  linearly constrained ones by introducing decision copies for interconnected agents.  This yields  linearly constrained nonconvex problems  with  local non-linear constraints.  
 The work \cite{sun2019two} argues that  the direct extension of ADMM to the reformulated problem is not applicable for the two necessary conditions \textbf{condition} a) and b) can not be satisfied simultaneously. 
 To bypass the challenge,   \cite{sun2019two} proposed  to introduce a block of slack variables working as the last block. To force the slack block  to \emph{zero}, this paper developed  a two-level  method where the inner-level uses classic ADMM  to  solve a relaxed problem with a penalty on the slack variables , and the outer-level gradually forces the slack variables  towards \emph{zero}.

 As can be perceived from the literature,  it is difficult (if not impossible) to develop a distributed  method with convergence guarantee for  \eqref{pp:problem template} due to the lack of a well-behaved last block satisfying \textbf{condition} a) and b). 
The work \cite{chatzipanagiotis2017convergence} provided  a solution with local convergence guarantee but can not handle the probable composite objective components $g$. 
 Though the idea of introducing slack variables in  \cite{sun2019two} can  provide a solution with global convergence guarantee but at the cost of  heavy iteration complexity caused by the two-level structure.  
 Despite these limitations, what we can learn from the literature  is that the behaviors of  dual variables is important to draw  the convergence of ADMM for nonconvex  problems.

 This paper focuses on developing a distributed  method for problem \eqref{pp:problem template} with theoretical convergence guarantee.  Our main contributions are 
%
%
%
\begin{itemize}
	\item We propose a proximal  ADMM  by revising the dual update procedure  of classic ADMM into a discounted manner. 
	 This leads to the boundness of  dual updates, which is critical to establish the convergence. 	
	\item We establish the global convergence of the method towards approximate stationary points by identifying  a proper  Lyapunov function which is sufficiently decreasing and lower bounded as required.
	\item We showcase the performance of the distributed method with a numerical example and a concrete application arising from smart buildings, which demonstrate  the method's effectiveness. 
\end{itemize}
The reminder of this paper is organized as follows.  In Section 2, we present the proximal  ADMM.  In Section 3, we study  the convergence of the method.  In Section 4, we showcase the method's performance with a numerical example and  smart building application.  In Section 5, we conclude this paper  and discuss the future work.

\section{Proximal ADMM}

\subsection{Notations}
Throughout the paper, we will visit the following notations. We use the bold alphabets $\x, \y, \mathbf{a}, \mathbf{b}, \mathbf{c}$ and  $\A, \A_{i}, \Q, \mathbf{M}$ to represent vectors and  matrices. We define $\mathbf{I}_n$ or $\mathbf{I}$ as identity matrices of $n \times n$ or suitable size.  We use the operator $:=$ to give  definitions.
We have $\R^{n}$  represent the $n$-dimensional real space and $(\x_i)_{i = 1}^{N}: = (\x_1^{\top}, \x_{2}^{\top}, \cdots, \x_{N}^{\top})^{\top} $ is the stack  of sub-vector $\x_{i} \in \R^{n_i}$. 
We refer to  $\Vert\!\cdot\!\Vert$ as Euclidean norm  without specification, i.e., $\Vert \x\Vert = \sqrt{{\textstyle \sum}_{i=1}^n x_i^2}$ for $\x \in \R^n$,  and $\langle \x, \y\rangle$ denote the dot product of vector $\x, \y \in \R^n$. We besides have $\Vert \x \Vert^{2}_{\A} = \x^{\top}\A \x$.
We use $\text{diag}(\A_{1}, \A_{2}, \cdots, \A_{N})$ to denote the diagonal matrix formed by the sub-matrices $\A_{1}, \A_{2}, \cdots, \A_{N}$.     We have the normal cone to  a convex set $\X \subseteq \R^n$ at $\starx$ defined by  $N_{\X}(\starx): = \{\nu \in \R^n \vert \langle \nu, \x - \starx \rangle \leq 0, \forall \x \in \X\}$. For $g: \R^{n} \rightarrow \R$ and $\x =(\x_i)_{i=1}^N \in \R^n$, we denote $\nabla_{i} g(\x) = \nabla_{\x_i} g(\x)$  as the partial differential of $g$ with respect to component $\x_i \in \R^{n_i}$.   We define  ${\rm dist}(\x, \X) = \min_{\y \in \X} \Vert \x - \y \Vert$  as the distance of vector $\x \in \R^n$ to the subset $\X \subseteq \R^n$. 
\subsection{Algorithm}

In this part, we introduce the proximal  ADMM for solving problem \eqref{pp:problem template} in a distributed manner. 
The proximal  ADMM is a type of  AL  methods that depend on the AL technique to relax  constraints and  employ the primal-dual scheme to update variables. 
By defining Lagrangian multipliers $\boldlambda \in \R^m$ for the coupled constraints \eqref{eq:1a}, we have the AL function for problem \eqref{pp:problem template}
\begin{align}
	\Lag_{\rho}(\x, \boldlambda) =& F(\x)+ \left\langle  \boldlambda,  \A\x -\bb \right\rangle \!+\! \frac{\rho}{2} \Vert \A \x -\bb  \Vert^2
\end{align}
where $F(\x) = g(\x) + \sum_{i =1}^N f_i(\x_i)$ and $\rho$ is the penalty parameter.   

Following the standard AL methods,   the proximal ADMM is  composed of \texttt{Primal update} and \texttt{Dual update} 
as shown in \textbf{Algorithm}  \ref{alg: perturbed_ADMM}. 
In \texttt{Primal update}, the primal variables $\x = (\x_i)_{i=1}^N$ are updated in a distributed manner via  \emph{Jacobian} decomposition.  Particularly, to handle the composite  objective component $g$, we  linearize the composite term  at each iteration $k$ by $g(\x^k) + \langle \nabla g(\x^k), \x - \x^{k} \rangle$ (the constant part $g(\x^k)$ is dropped).  Note that the local objective terms $f_i$ can also be linearized similarly if necessary and the proof of this paper still applies.
To favor computation efficiency and scaling properties,  we adopt the \emph{Jacobian} scheme and empower the agents to update their decision components  in parallel at each iteration with the preceding information from their  interconnected agents.  
Particularly, to enhance  convergence,  a  proximal term $\Vert \x_i - \x_i^{k + 1} \Vert^2$ is imposed on the local objective of each agent (Step 3). This has been  used in many \emph{Jacobian} ADMM  both in convex \cite{li2020distributed, deng2017parallel, chang2014multi}  and  nonconvex  \cite{liu2019linearized, lu2021linearized} settings.   
Note that the subproblems \eqref{eq:x-update} are either  convex or nonconvex optimization over the local constraints $\X_i$, depending on $f_i$.  There are many first-order solvers to solve those subproblems, such as  the projected gradient method \cite{jain2017non}   and the proximal  gradient method \cite{li2015accelerated}. This paper focuses on developing a general distribute framework for solving problem \eqref{pp:problem template} and will not discuss the subproblems in detail.
The major difference of the proximal ADMM from  the existing distributed AL methods is that we have modified  the \texttt{Dual update} by  imposing a  discounting factor $(1-\tau)$ ($\tau \in [0, 1)$) (Step 4).
 The idea and motivation behind are  to update the dual variables by the constraints residual in a discounted manner so as to bound the dual variables in the iterative process, which has been identified as critical to draw  theoretical convergence. 
  In this setting,  the dual variables are the \emph{discounted} running sum of the constraints residual, i.e.,  

 {\small 
 \begin{equation} \label{eq:lambda}
	\setlength{\abovedisplayskip}{-2pt}
	\setlength{\belowdisplayskip}{2pt}
\begin{split}
\boldlambda^{k +1} & = (1  -\tau) \boldlambda^k + \rho(\A \x^{k + 1}- \bb)\\ 
& = (1 - \tau)^2 \boldlambda^{k - 1 } + (1 - \tau) \rho (\A \x^k - \bb)\\
& \quad \quad  \quad \quad \quad \quad \quad + \rho (\A \x^{k + 1} - \bb)  \\
& \cdots \\
   & = (1 -\tau)^{k+1} \boldlambda^0 + \smallsum_{\ell=0}^{k} (1-\tau)^{k -\ell} \rho(\A\x^{\ell + 1} - \bb).
\end{split}
 \end{equation}
}
This differs from classic ADMM where  the dual variables are the running sum of the constraints residual, i.e., 
 {\small 
		\setlength{\abovedisplayskip}{-1pt}
		\setlength{\belowdisplayskip}{2pt}
\begin{align*}
    \boldlambda^{k + 1} & = \boldlambda^{k} + \rho (\A \x^{k + 1} - \bb) \\
                    & =  \boldlambda^{k - 1} + \rho (\A \x^{k } - \bb) + \rho (\A \x^{k + 1} - \bb) \\
                    & \cdots \\
                    & = \boldlambda^{0} + \smallsum_{\ell = 0}^{k} \rho (\A \x^{\ell + 1} - \bb).
\end{align*}}
From this perspective,  classic ADMM can be viewed as a special case of the proximal ADMM with $\tau = 0$.    
In the proximal ADMM, the \texttt{Primal update} and \texttt{Dual update} are alternated until the stopping criterion 
{	\setlength{\abovedisplayskip}{3pt}
	\setlength{\belowdisplayskip}{3pt}
\begin{align} \label{eq:stopping_criterion}
\Vert T_c^{k+1} - T_c^k \Vert \leq \epsilon
\end{align}}
is reached, where $T_c^k$ is the Lyapunov function to be discussed later. The parameter  $\epsilon$ is a user-defined  positive threshold.

\begin{algorithm}[h]
	\caption{Proximal ADMM for problem \eqref{pp:problem template} }
	\label{alg: perturbed_ADMM} 
	\begin{algorithmic}[1] 
		\State  \textbf{Initialize:} $\mathbf{x}^0$, $\mathbf{\boldlambda}^0$ and $\rho>0$, $\tau \in [0, 1)$,   and set $k \rightarrow 0$.
		\State \textbf{Repeat:} 
		\State ~~\texttt{Primal update}:
		{\small 
			\begin{equation}
				\begin{split}
					&\label{eq:x-update} \x_i^{k+1}\!= \!\arg\!\min_{\x_i \in \X_i} \!\!\left\{  
					\begin{array}{l}
						\langle \nabla_i g(\x^{k}), \x_{i} - \x_{i}^{k} \rangle \\
						\!+ f_{i}(\x_{i}) 
						+ \langle \boldlambda^{k}, \A_{i} \x_{i}^{k}\rangle \\
						\!+ \rho/2\Vert \A_{i} \x_{i} \!+\! \!{\textstyle\sum}_{j \neq i} \A_{j} \x_{j}^{k} \!-\! \bb \Vert^{2} \!\!\\
						\!+ \beta/2 \Vert \x_{i} - \x_{i}^{k} \Vert^{2}_{\B_{i}}
					\end{array}
					\right\} \\
				\end{split}
		\end{equation} }
		\State~~ \texttt{Dual update}:
		{\small 
			\begin{align}
				\label{eq:dual-update}\boldlambda^{k+1} & \!=\! (1-\tau)\boldlambda^k \!+\! \rho \left(   \A \x^{k+1} \!-\!\mathbf{b}  \right) 
		\end{align}}
		\State Until the stopping criterion \eqref{eq:stopping_criterion}  is reached. 
	\end{algorithmic}
\end{algorithm}

\section{Convergence Analysis}

Before establishing the convergence of \textbf{Algorithm}  \ref{alg: perturbed_ADMM}, we first clarify the main assumptions. 

\subsection{Main assumptions}
\begin{itemize}
\item[(A1)] Function $f\!:\!\R^{n} \rightarrow \R$ and $g: \R^{n} \rightarrow \R$ have $\Lips$ continuous gradient (i.e., $\Lips$ differentiable) with modulus $L_{f}$ and $L_{g}$ over the  set  $\X = \X_1 \times \X_2 \times \cdots \times \X_N$, i.e., \cite{guo2017convergence}
\begin{align*}
&\Vert \nabla f(\x) -  \nabla f(\y) \Vert   \leq L_{f} \Vert \x - \y \Vert,  ~\forall \x, \y \in \X. \\
   & \Vert \nabla g(\x) - \nabla g(\y) \Vert \leq L_{g} \Vert \x - \y \Vert,  ~\forall \x, \y \in \X. \\
\end{align*} 
	\item[(A2)] Function $\!f\!:\R^{n} \!\!\rightarrow\! \R$ and $g\!: \R^{n} \!\rightarrow\!\! \R$ are lower bounded over the set $\X \!=\! \X_1 \!\times\! \X_2 \times\! \cdots \!\times\! \X_N$, i.e., 
	\begin{align*}
	f(\x) & > - \infty, ~~\forall \x\in \X. \\
	g(\x) & > -\infty,~~ \forall \x \in \X. 
	\end{align*}
\end{itemize}





\subsection{Main results}
As discussed, there are two key steps to draw convergence for a distributed AL method in nonconvex settings: 1)  identifying  a so-called sufficiently decreasing Lyapunov function; and  2) establishing the lower boundness property of the Lyapunov function. To achieve the objective, we first draw the  following two propositions.

%
\begin{proposition} \label{prop:prop-x}
	For the sequences $\{\x^k\}_{k \in \K}$ and $\{\boldlambda^k\}_{k \in \K}$ generated by \textbf{Algorithm}  \ref{alg: perturbed_ADMM},  we have
\begin{equation}
	\begin{split}
	& \frac{1-2\tau^2}{2\rho}   \left\Vert \boldlambda^{k+1} -\boldlambda^k \right\Vert^2   + \frac{1}{2}\Vert  \x^{k+1}-\x^k\Vert^2_{\Q} \\
& \quad + \frac{L_g}{2}\Vert \x^{k +1} - \x^{k}\Vert^{2} + \frac{1}{2} \Vert \w^k \Vert^2_{\Q}\notag \\
&  \leq  \frac{1-2\tau^2}{2\rho}   \left \Vert  \boldlambda^k-\boldlambda^{k-1} \right\Vert^2 + \frac{1}{2}\Vert  \x^k-\x^{k-1}\Vert^2_{\Q} \\
&\quad + \frac{L_g}{2}\Vert \x^{k} - \x^{k - 1}\Vert^{2} + \rho_F\left\Vert \x^{k+1} -\x^k \right\Vert^2 \\
& \quad- {\tau(1 + \tau)}/{\rho} \Vert \boldlambda^{k+1}-\boldlambda^k\Vert^2.
	\end{split}
\end{equation}
where we have the iterations $\K:=\{1, 2, \cdots, K\}$  and 
\vspace{3mm}
{\small 
\setlength{\abovedisplayskip}{-30pt}
	\setlength{\belowdisplayskip}{3pt}
\begin{align*}  \\
     \w^k:   & = (\x^{k+1}-\x^k) - (\x^k-\x^{k-1}) \\
	G_{\A}:  & = \text{diag}\left(\A_1^\top \A_1, \cdots, 	\A_N^\top \A_N  \right) \\
	G_{\B}:  & = \text{diag} \left( \B_1^\top \B_1 , \cdots, \B_N^\top \B_N \right) \\
    \Q:      & =\rho G_{\A} + \beta G_{\B} - \rho\A^\top\A \\
    \rho_F: & = L_{f} + L_{g}.
\end{align*}}
\end{proposition}
\textbf{Proof of Prop. \ref{prop:prop-x}}: 
We defer the proof to \textbf{Appendix} A.

Let $\Lag^{+}_{\rho}(\x, \boldlambda) := \Lag_{\rho}(\x, \boldlambda) - \frac{\tau}{2\rho}\Vert \boldlambda\Vert^2$  be the regularized AL function. We have the subsequent proposition to quantify the change of   regularized AL function over the successive iterations. 
\begin{proposition} \label{prop:prop-L}
	For the sequences $\{\x^k\}_{k \in \K}$ and $\{\boldlambda^k\}_{k \in \K}$ generated by \textbf{Algorithm}  \ref{alg: perturbed_ADMM},  we have
\begin{equation*}
	\setlength{\abovedisplayskip}{4pt}
   \setlength{\belowdisplayskip}{1pt}
	\begin{split}
		& \Lag^{+}_{\rho}(\x^{k+1}, \boldlambda^{k + 1})  -  \Lag^{+}_{\rho}(\x^k, \boldlambda^k) \\
		&  \leq  - \Vert \x^{k+1} - \x^k \Vert_{\Q}^2 + \frac{\rho_F}{2} \Vert \x^{k+1} -\x^{k} \Vert^2\!\\
		& \quad  -\frac{\rho}{2} \Vert \A(\x^{k+1} - \x^k)\Vert^2  + \frac{2-\tau}{2\rho} \Vert \boldlambda^{k+1} - \boldlambda^k\Vert^2.
	\end{split}
\end{equation*}
\end{proposition}
\textbf{Proof of Prop. \ref{prop:prop-L}:}  We defer the proof to \textbf{Appendix} B.

In the literature, the AL function is often  used  as the Lyapunov function if the   sufficiently decreasing property can be established (see \cite{wang2019global, yang2017alternating, guo2017convergence, li2015global} for examples).   
However, this is not the case for \textbf{Algorithm} \ref{alg: perturbed_ADMM}. 
From \textbf{Prop.} \ref{prop:prop-L}, we note that the sufficiently decreasing property of the (regularized) AL function can be established if and only if the dual updates $\Vert \boldlambda^{k + 1} - \boldlambda^{k} \Vert^2$ 
can be bounded by the primal updates $\Vert \x^{k+1} - \x^k \Vert^2$ (see the definition \eqref{eq:sufficiently-decreasing}). 
This is difficult (if not impossible)  due to the lack of a well-behaved last block (i.e., unconstrained and $\Lips$ differentiable) as discussed. 

However, by combing \textbf{Prop.} \ref{prop:prop-x} and \textbf{Prop.} \ref{prop:prop-L}, we indeed can identify a sufficiently decreasing Lyapunov function.  Specifically,  from \textbf{Prop.} \ref{prop:prop-L}, we have the (regularized) AL function  $\Lag^{+}_{\rho}(\x^{k+1}, \boldlambda^{k + 1})$ is ascending in $\Vert \boldlambda^{k+1} - \boldlambda^k\Vert^2$  and descending in  $\Vert \x^{k+1} - \x^k \Vert^2$.  This is exactly  opposite to the descending and ascending  properties of the term  $ \frac{1-2\tau^2}{2\rho}  \Vert \boldlambda^{k+1} -\boldlambda^k \Vert^2 + \frac{1}{2}\Vert  \x^{k+1}-\x^k\Vert^2_{\Q}$ stated in \textbf{Prop.} \ref{prop:prop-x}. 
Note that this is attributed to  the imposed discounted factor $\tau > 0$, otherwise the term $\tau(1-\tau)/\rho\Vert \boldlambda^{k + 1} - \boldlambda^k \Vert^2$ in \textbf{Prop.} \ref{prop:prop-x} would be \emph{zero}.
We  therefore build  the  Lyapunov function as 
{\setlength{\abovedisplayskip}{4pt}
	\setlength{\belowdisplayskip}{3pt}
\begin{align} \label{eq:T-def}
	& T_c(\x^{k+1}, \boldlambda^{k+1}; \x^k, \boldlambda^k)  = \Lag^{+}_{\rho}(\x^{k+1}, \boldlambda^{k + 1})\notag\\
	&  \quad ~ + c \bigg(  \frac{1-2\tau^2}{2\rho}  \Vert \boldlambda^{k+1} -\boldlambda^k \Vert^2 + \frac{1}{2}\Vert  \x^{k+1}-\x^k\Vert^2_{\Q} \notag \\
	& \quad ~+ \frac{L_g}{2}\Vert \x^{k} - \x^{k-1} \Vert^{2} \bigg)  
\end{align}}
where $c$ is a constant parameter to be determined  for ensuring  the sufficiently decreasing and lower boundness property of  the Lyapunov function.


Let  $T_c^{k+1}: = T_c(\x^{k+1}, \boldlambda^{k+1}; \x^k, \boldlambda^k)$ be the Lyapunov function at iteration $k$, we have the following proposition regarding the sufficiently decreasing property. 
\begin{proposition}\label{prop:sufficient-decrease}
	For the sequences $\{\x^k\}_{k \in \K}$ and $\{\boldlambda^k\}_{k \in \K}$ generated by \textbf{Algorithm}  \ref{alg: perturbed_ADMM}, we have  
	\begin{equation*} 
		\begin{split}
			& T_c^{k+1}\!- \!T_c^k \leq -a_{\x} \Vert \x^{k+1} \!-\! \x^k\Vert^2 \!-\! a_{\boldlambda} \Vert \boldlambda^{k+1} \!-\! \boldlambda^k \Vert^2 \!-\!\frac{c}{2}\Vert \w^k \Vert^2
		\end{split}
	\end{equation*}
	where we have $\rho_F = L_{f} + L_{g}$ and
	\begin{equation*} \label{eq:definition}
		\begin{split}
			& a_{\x}: = \frac{ 2\rho G_{\A} + 2\beta G_{\B} -\rho \A^\top\A  - (2c + 1)\rho_F \mathbf{I}_N }{2}  \\
			& a_{\boldlambda}: = \frac{2c\tau(1 + \tau) - (2-\tau )}{2\rho}.
		\end{split}
	\end{equation*}
\end{proposition}
\textbf{Proof of Prop.  \ref{prop:sufficient-decrease}}: Based on \textbf{Prop.}  \ref{prop:prop-x} and \textbf{Prop.} \ref{prop:prop-L}, we have
\begin{align*}
	& T_c^{k+1} - T_c^k  =  - \Vert \x^{k+1} - \x^k \Vert_{\Q}^2 + \frac{\rho_F}{2} \Vert \x^{k+1} -\x^{k} \Vert^2   \notag \\
	& \quad   -\frac{\rho}{2} \Vert \A(\x^{k+1} - \x^k)\Vert^2 + \frac{2-\tau}{2\rho} \Vert \boldlambda^{k+1} - \boldlambda^k\Vert^2  \notag\\
	& \quad + c \Big(\rho_F \left\Vert \x^{k+1} -\x^k \right\Vert^2 - {\tau(1 + \tau)}/{\rho} \Vert \boldlambda^{k+1}-\boldlambda^k\Vert^2  \notag \\
	& \quad- {1}/{2} \Vert \w^k \Vert^2_{\Q} \Big)   \notag\\
	& \leq -a_{\x}\Vert \x^{k + 1} - \x^k  \Vert^2  - a_{\boldlambda} \Vert  \boldlambda^{k +1} - \boldlambda^k \Vert^2 - \frac{c}{2}\Vert \w^k \Vert^2_{\Q}
\end{align*}
where the inequality is directly derived by rearranging the terms. We therefore close  the proof. 

\begin{remark}
	\textbf{Prop.} \ref{prop:sufficient-decrease} implies that we would  have the sufficiently decreasing property hold by  the constructed Lyapunov function  $T_c^{k}$  if we have $a_{\x} > 0$, $a_{\boldlambda} > 0$,  $c \geq 0$ and $\Q \geq 0$.  Actually, this can be achieved by setting the tuple  ($\tau$, $\rho$, $\beta$, $\B_i$, $c$)  properly for \textbf{Algorithm} \ref{alg: perturbed_ADMM}, which will be discussed shortly. 
\end{remark}


As discussed,  another  key step to draw the convergence is to establish  the lower boundness property of the Lyapunov  function.  To this end,  we first prove the lower boundness property of  Lagrangian multipliers resulting from the discounted dual  update scheme.  

\begin{proposition} \label{prop:lambda-bounded}  Let $\Delta^k: = \Vert  \A \x^k -\bb \Vert $ be the constraints residual  at iteration $k$,  $\Delta^{\max}: = \max_{\x \in \X} \Vert \A\x -\bb \Vert$  denote the maximal  constraints residual over the closed  feasible set $\X$, and \textbf{Algorithm} \ref{alg: perturbed_ADMM}  start  with any given initial dual variable $\boldlambda^0$,  we have $\Vert \boldlambda^k \Vert$ is bounded, i.e., 
	\begin{equation} \label{eq:lambda-bound}
		\setlength{\abovedisplayskip}{4pt}
		\setlength{\belowdisplayskip}{3pt}
		\begin{split}
			& \Vert \boldlambda^k \Vert \leq \Vert \boldlambda^0 \Vert + \tau^{-1}\rho\Delta^{\max}  \\
			\text{or}~ ~& \Vert \boldlambda^k\Vert^2 \leq 2\Vert \boldlambda^0 \Vert^2 + 2\tau^{-2}\rho^2 (\Delta^{\max})^2.
		\end{split}
	\end{equation}
\end{proposition}
\textbf{Proof of \textbf{Prop.} \ref{prop:lambda-bounded}:}
Recall  the dual update scheme in  \eqref{eq:lambda}, we have

	
{\setlength{\abovedisplayskip}{-3pt}
	\setlength{\belowdisplayskip}{3pt} 
		\begin{align*}
			\Vert \boldlambda^{k} \Vert & = \Vert (1-\tau)^{k + 1} \boldlambda^0 +  \smallsum_{\ell=0}^k \rho (1-\tau)^{k-\ell} \Delta^{\ell + 1} \Vert \\
& \leq \Vert (1-\tau)^{k + 1} \boldlambda^0 \Vert +\smallsum_{\ell = 0}^k \Vert \rho (1-\tau)^{k - \ell} \Delta^{\ell + 1}  \Vert \\
& \leq \Vert (1-\tau)^{k + 1} \boldlambda^0 \Vert  + \rho\Delta^{\max} \frac{1-(1-\tau)^k}{\tau} \\ 
			& \leq \Vert \boldlambda^0 \Vert + \tau^{-1}\rho\Delta^{\max} 
	\end{align*}}
where the first inequality   is by the triangle inequality of norm,  the second inequality  infers from $\Delta^k \leq \Delta^{\max}, \forall k$, and  the last inequality holds because of $\tau \in (0, 1)$. 
	
	Further based on $\Vert \mathbf{a} + \mathbf{b}\Vert ^2 \leq 2\Vert \mathbf{a} \Vert^2 + 2\Vert \mathbf{b} \Vert^2$, we directly  have $ \Vert \boldlambda^{k}\Vert^{2} \leq 2\Vert \boldlambda^{0}\Vert^{2} + 2\tau^{-2}\rho^2 (\Delta^{\max})^2$. We therefore complete the proof. 
	
	Based on \textbf{Prop.} \ref{prop:lambda-bounded}, we are able to establish  the lower boundness property  of  Lyapunov function  as below. 
	\begin{proposition} \label{prop:lower-bounded}
		For the sequences $\{\x^k\}_{k \in \K}$ and $\{\boldlambda^k\}_{k \in \K}$ generated by \textbf{Algorithm}  \ref{alg: perturbed_ADMM}, we have 
		{\setlength{\abovedisplayskip}{3pt}
			\setlength{\belowdisplayskip}{3pt} 
		\begin{align} \label{eq:lower-bounded}
			T_c^{k+1} > - \infty, \forall k \in \K.
		\end{align}}
	\end{proposition}
	\textbf{Proof of Prop. \ref{prop:lower-bounded}:} By examining the terms of $T_c^{k+1}$  in \eqref{eq:T-def},  we only require  to  establish the lower boundness property of  $\Lag_{\rho}^{+}(\x^{k+1}, \boldlambda^{k+1}) = \Lag_{\rho}(\x^{k+1}, \boldlambda^{k+1}) - \frac{\tau}{2\rho} \Vert \boldlambda^{k+1} \Vert^2$ for the other terms are all non-negative.  
	Based on \textbf{Prop.} \ref{prop:lambda-bounded}, we directly have 
	 $- \frac{\tau}{2\rho}\Vert \boldlambda^{k+1} \Vert^2 $ lower bounded since $\Vert \boldlambda^{k+1} \Vert^2 $ is upper bounded. 
	We therefore only need to prove that
	$ \Lag_{\rho}(\x^{k+1}, \boldlambda^{k+1}) = f(\x^{k+1})  + \langle \boldlambda^{k+1}, \A\x^{k+1}-b  \rangle 
	+{\rho}/{2}\left\Vert  \A\x^{k+1}-\bb\right\Vert^2 $
	is lower bounded. Note that we have $f(\x^{k+1}) > -\infty$ over the compact set $\X$ (see (A2)) and the quadratic term non-negative. This infers we only need to prove the lower boundness for the second term $ \langle \boldlambda^{k+1}, \A\x^{k+1}-b  \rangle $.  Based on the dual update \eqref{eq:dual-update}, we have 
	\begin{align}\label{eq:second-term_L}
		& \langle \boldlambda^{k+1}, \A\x^{k+1}- \bb\rangle = \left\langle \boldlambda^{k+1}, \frac{\boldlambda^{k+1}-(1-\tau)\boldlambda^k}{\rho}\right\rangle \notag \\
		& =  \left\langle  \boldlambda^{k+1}, \frac{1-\tau}{\rho}(\boldlambda^{k+1}-\boldlambda^k) + \frac{\tau}{\rho} \boldlambda^{k+1}\right\rangle  \\
		& = \frac{\tau}{\rho} \Vert \boldlambda^{k+1}\Vert^2 + \frac{1-\tau}{\rho} \left\langle \boldlambda^{k+1}, \boldlambda^{k+1}-\boldlambda^k\right\rangle  \notag  \\
		& = \frac{\tau}{\rho} \Vert \boldlambda^{k+1} \Vert^2 \!+\! \frac{1-\tau}{2\rho}\big( \Vert \boldlambda^{k+1}\!-\!\boldlambda^k\Vert^2 \!+\! \Vert \boldlambda^{k+1}\Vert^2 \!-\! \Vert \boldlambda^k \Vert^2\big) \notag
	\end{align}
	Since we have  $\Vert \boldlambda^{k} \Vert^{2}$ is upper bounded (see \textbf{Prop.} \ref{prop:lambda-bounded}), we therefore have  $\Lag_{\rho}(\x^{k+1}, \boldlambda^{k+1}) $ lower bounded for  the other terms of \eqref{eq:second-term_L}  are all non-negative. We thus complete the proof.

	To present the main results regarding the convergence of \textbf{Algorithm} \ref{alg: perturbed_ADMM}, we first give the definition on \textbf{Approximate stationary solution}.  
	\begin{definition}
		\textbf{(Approximate stationary solution)} For any given $\epsilon$, we say a tuple $(\x^{*}, \boldlambda^{*})$ is an $\epsilon$-stationary solution of problem \eqref{pp:problem template}, if we have 
		\begin{equation*}
			\begin{split}
				& \text{\emph{dist}}\big(\nabla F(\starx)+ \A^\top \boldlambda^{*} + N_{\X}(\starx), \mathbf{0}\big)+ \Vert \A \starx- \bb\Vert  \leq \epsilon. \\
			\end{split}	
		\end{equation*}
	\end{definition}
	where $\nabla F(\x^{*}) = \nabla f(\x^{*}) + \nabla g(\x^{*})$. 
	
	In terms of the convergence of \textbf{Algorithm} \ref{alg: perturbed_ADMM} for problem \eqref{pp:problem template}, we have the following main results. 
	\begin{theorem} \label{thm:theorem1}
		For \textbf{Algorithm} \ref{alg: perturbed_ADMM} with the tuple  ($\tau$, $\rho$, $\beta$, $\B_i$, $c$) selected by 
		{\setlength{\abovedisplayskip}{3pt}
			\setlength{\belowdisplayskip}{3pt}
\begin{align}
			 & \tau:   \tau \in (0, 1)  \notag\\ 
			  & c:    c > \frac{2 - \tau}{2 \tau(1 + \tau)}\tag{C1}\\
				& (\rho, \beta, \B_i):\!
				  \begin{cases} 
			       & \!\!\!\!\! 2\rho G_{\A} \!+\! 2 \beta G_{\B} \!-\! \rho \A^\top \A \geq (2 c +1) \rho_F \mathbf{I}_N \notag \\
				 &\!\!\! \!\! \Q: = \!\rho G_{\A} + \beta G_{\B} \!-\! \rho \A^\top\A \geq 0   \notag   
				\end{cases}
\end{align}
		}
		\begin{itemize}
			\item[(a)] The generated sequence $\{\x^k\}_{k \in \K}$ and $\{\boldlambda^k\}_{k \in \K}$ are bounded and convergent, i.e.,
			{\setlength{\abovedisplayskip}{5pt}
				\setlength{\belowdisplayskip}{-10pt}
			\begin{align*}
				\boldlambda^{k+1} - \boldlambda^{k} \rightarrow 0, ~~\x^{k+1} - \x^k \rightarrow 0. \\
			\end{align*}}
			\item[(b)]   Suppose we have the limit tuple $(\x^{*}, \boldlambda^{*})$,   then $(\x^{*}, \hat{\boldlambda}^{*})$  with $\hat{\boldlambda}^{*} = (1 + \tau \boldlambda^{*})$ is  $\tau \rho^{-1} \Vert \boldlambda^{*} \Vert$-stationary solution of problem \eqref{pp:problem template}. 
		\end{itemize}
	\end{theorem}
	
	\textbf{Proof of Theorem \ref{thm:theorem1}}: (a) Recall  \textbf{Prop.} \ref{prop:sufficient-decrease}, we have 
	{\small 
		\begin{equation*}
			\setlength{\abovedisplayskip}{-2pt}
			\setlength{\belowdisplayskip}{1pt}
			\begin{split}
				& \sum_{k=1}^{K} \big( T_c^k - T_c^{k+1} \big)  \geq  a_{\x} \sum_{k=1}^{K} \Vert \x^{k+1} \!-\! \x^k\Vert^2  \\
				& \quad \quad \quad \quad \quad \quad+ a_{\boldlambda} \sum_{k=1}^{K}\Vert \boldlambda^{k+1} \!-\! \boldlambda^k \Vert^2+ \frac{c}{2} \sum_{k=1}^{K} \Vert \w^k \Vert^2 \\
							\end{split}	
						\end{equation*}}
		By assuming $K \rightarrow \infty$, we have
						\begin{equation*}
							\setlength{\abovedisplayskip}{-2pt}
							\setlength{\belowdisplayskip}{1pt}
							\begin{split}
				&T_c^1 - \lim_{K \rightarrow \infty}T_c^{k+1} \geq a_{\x} \sum_{k=1}^{\infty}\Vert \x^{k+1} \!-\! \x^k\Vert^2 \\
				&  \quad \quad \quad \quad \quad \quad+ a_{\boldlambda} \sum_{k=1}^{\infty}\Vert \boldlambda^{k+1} \!-\! \boldlambda^k \Vert^2  + \frac{c}{2} \sum_{k=1}^{\infty} \Vert \w^k \Vert^2 \\
				\end{split}	
				\end{equation*}
				Since we have $T_c^{k+1} > -\infty$ (see \textbf{Prop.} \ref{prop:lower-bounded}), we thus have 
				\vspace{-3mm}
				{\small
								\begin{equation*}
									\begin{split}
		& \infty \!\geq\! a_{\x} \sum_{k=1}^{\infty}\Vert \x^{k+1} \!-\! \x^k\Vert^2 \!+\! a_{\boldlambda} \sum_{k=1}^{\infty}\Vert \boldlambda^{k+1} \!-\! \boldlambda^k \Vert^2 \!+\!  \frac{c}{2} \sum_{k=1}^{\infty} \Vert \w^k \Vert^2.  \\
			\end{split}	
		\end{equation*}}
	We therefore conclude
	{\small 
		\begin{equation*}
			\setlength{\abovedisplayskip}{-2pt}
			\setlength{\belowdisplayskip}{1pt}
			\begin{split}
				& \Vert \x^{k+1} - \x^k \Vert \rightarrow 0, ~~~ \Vert \boldlambda^{k+1} - \boldlambda^k \Vert \rightarrow 0, \\
				& \Vert \w^k \Vert =\Vert (\x^{k+1}-\x^k) - (\x^k-\x^{k-1}) \Vert \rightarrow 0. 
			\end{split}	
		\end{equation*}
	}

	(b) According to (a), we have the sequences $\{\x^k\}_{k \in \K}$ and $\{\boldlambda^k\}_{k \in \K}$   converge to some limit tuple $(\x^{*}, \boldlambda^{*})$, i.e., if $k \rightarrow \infty$, we have $\x^{k+1} \rightarrow \starx, \boldlambda^{k+1} \rightarrow \boldlambda^{*}$ and $\x^{k+1} \rightarrow \x^k$ and $\boldlambda^{k+1} \rightarrow \boldlambda^{k}$. 
	
	Based on the dual update procedure \eqref{eq:dual-update},  we have the stationary tuple $(\starx, \boldlambda^{*})$ satisfy 
		\begin{equation}
		\setlength{\abovedisplayskip}{5pt}
			\setlength{\belowdisplayskip}{3pt}
			\begin{split} \label{eq:inequality6}
				&\A \starx - \bb = \tau \rho^{-1} \boldlambda^{*}. \\
			\end{split} 
		\end{equation}
Since we have  $\hat{\boldlambda}^k = \boldlambda^k + \rho(\A \x^k - \bb)$,  we thus have $\hat{\boldlambda}^k \rightarrow (1 + \tau) \boldlambda^{*}$. Let $\hat{\boldlambda}^{*} = (1 + \tau) \boldlambda^{*}$, we have $\hat{\boldlambda}^k \rightarrow \hat{\boldlambda}^{*}$. 
	

Recall the first-order optimality condition \eqref{eq:first-order-k} and assume $k \rightarrow \infty$ that  the stationary point ($\starx, \boldlambda^{*}$) is reached,  we would have 
	\begin{align*}
\langle \nabla f(\starx) \!+\! \nabla g(\starx) \!+\! \A^\top\hat{\boldlambda}^{*}, \starx - \x \rangle  \leq 0, \forall \x \in \X. 
	\end{align*}	
	This implies that 
	\begin{align*}
 \bm{0} \in \nabla f(\starx)  \!+\! \nabla g(\starx)  \!+\! \A^\top \hat{\boldlambda}^{*} + N_{\X}(\starx). 
	\end{align*}
We further have 
	\begin{equation}\label{eq:approximate}
		\begin{split}
			& \text{dist}\big( \nabla  f(\starx)  \!+\! \nabla g(\starx) \!+\! \A^\top  \hat{\boldlambda}^{*} + N_{\X}(\starx) , 0\big)=0\\
		\end{split}
	\end{equation}
	By combing \eqref{eq:inequality6} and \eqref{eq:approximate}, we therefore conclude
	\begin{align*}
& \text{dist}\big( \nabla  f(\starx) + \nabla g(\starx) + \A^\top \hat{\boldlambda}^{*} + N_{\X}(\starx) , 0\big) \\
&+ \Vert \A \starx - \bb  \Vert  \leq  \tau \rho^{-1} \Vert \boldlambda^{*} \Vert, 
	\end{align*} 
which closes the proof.

	From \textbf{Theorem} \ref{thm:theorem1}, we note that if  the convergent $\boldlambda^{*}$ does not depend on $\tau$ and $\rho$, we could decrease $\tau$  or increase $\rho$ to achieve any sub-optimality.  If that is not the case, we give the following corollary to show that this still can be  achieved  by 
	properly setting the initial point and parameters.

\begin{corollary}\label{cor:corollary}
	For any given $\epsilon > 0$,   if  \textbf{Algorithm} \ref{alg: perturbed_ADMM} starts with  $\boldlambda^0 = 0$ and $\A \x^0 = \bb$,  $\tau \in (0, 1)$, and  the penalty parameter $\rho$ is selected that  
	\begin{align*} 
		&\rho \geq ~ \epsilon^{-2} \tau  \big( 4 + c(1 - 2\tau^2) + c/2\big) d_F+ \epsilon^{-2} {c L_g}/{2}\Vert \x^0 \Vert^2 \\
		& \quad \quad + \epsilon^{-2}\tau {c L_g}/{2}\Vert \x^0 \Vert^2 +  \epsilon^{-2}\tau  {c \rho_F}/{4}~d_{\x},
	\end{align*}
 we have the limit tuples $(\x^{*}, \boldlambda^{*})$ and $(\x^{*}, \hat{\boldlambda}^{*})$  with $\hat{\boldlambda}^{*} = (1 + \tau \boldlambda^{*})$ is $\epsilon$-stationary solution of problem \eqref{pp:problem template}.
	where we have  $d_F= \max_{\x \in \X} f(\x) + g(\x)$, $d_{\x} = \max_{\x, \y \in \X} \Vert \x - \y\Vert^2$, and we assume $f(\x) \geq 0$, $g(\x) \geq 0$ without losing any generality. 
	\end{corollary}

\textbf{Proof of Corollary \ref{cor:corollary}}:  We only give the sketch of the proof and defer the details to \textbf{Appendix} C. The proof is structured by two parts which include: i) proving $\Vert \boldlambda^{*} \Vert^2 \leq  \rho \tau^{-1} T_c^0$,  and ii) proving   $T_c^0 \leq  \big( 4 + c(1 - 2\tau^2) + c/2\big) d_F+ {c L_g}/{2}\Vert \x^0 \Vert^2 +  {c \rho_F}/{4}~d_{\x}$. Based on i) and ii), we  have 
$\tau^2 \rho^{-2} \Vert \boldlambda^{*} \Vert^2  \leq \epsilon^2$. We then directly draw the conclusion 
based on \textbf{Theorem} \ref{thm:theorem1}.

\section{Numerical Experiments}
\subsection{A numerical example}

We first consider a numerical example with $N=2$ agents given by 
\begin{align}
	\label{pp:numeric_example}\min_{x_1, x_2} ~~& 0.1x_1^3 + 0.1x_2^3 + 0.1 x_1 x_2  \tag{$\mathbf{P1}$}  \\
	\text{s.t.} ~ &~ x_1 +  x_2 = 1 \notag\\
	& -1 \leq x_1 \leq 1 \notag\\
   &  -1  \leq x_2 \leq 1 \notag
\end{align}

For this example, we have  $f_1(x_1) =0.1 x_1^3$, $f_2(x_2) =0.1 x_2^3$,  and $g(x_1, x_2) = 0.1 x_1 x_2$. The $\Lips$ continuous gradient modulus for $f$  and $g$ are $L_f = 0.6$ and $L_g = 0.2$. Besides, we have  $\A_1 = 1$, $\A_2 = 1$, 
$\A = (1 ~ 1)$.  The stationary point of the problem is $x^{\star}_1= 0.5, x^{\star}_2 = 0.5$.

To our best knowledge, there is no distributed solution methods  for solving problem  \eqref{pp:numeric_example} with theoretical convergence guarantee.   In the following, we apply the proximal ADMM to solve this problem and verify the solution quality.  We consider four different  parameter settings  for  \textbf{Algorithm} \ref{alg: perturbed_ADMM}: 
\begin{itemize}
\item [] S1) $\tau = 0.1$,  $\rho = 10$, $\beta=10$,  $c = 8.7$
\item [] S2) $\tau = 0.1$, $\rho = 20$, $\beta = 20$,  $c= 8.7$
\item [] S3) $\tau =0.05$, $\rho = 5$, $\beta = 16$, $c = 18.6$
\item [] S4) $\tau = 0.05$, $\rho =10$, $\beta = 16$, $c = 18.6$
\end{itemize}
The other parameters are set as   $\B_1 = \B_2 =1$, $\tau = 0.1$,   $x_1^0 = 0.2$,   $x_2^0 = 0.8 $, $\lambda^0 = 0$ and  kept the same for S1-S4.  
Note that  we have $\tau/\rho = 0.01$ for S1/S3 and $\tau/\rho = 0.005$ for S2/S4.  We make such settings for comparisons as  we have the suboptimality of the method related to 
the ratio  $\tau/\rho$ as stated in \textbf{Theorem} \ref{thm:theorem1}. We therefore study how the ratio $\tau/\rho$ will affect the convergence rate and the solution quality of the method.

\begin{figure}[h]
 	\centerline{\includegraphics[width=3.0 in]{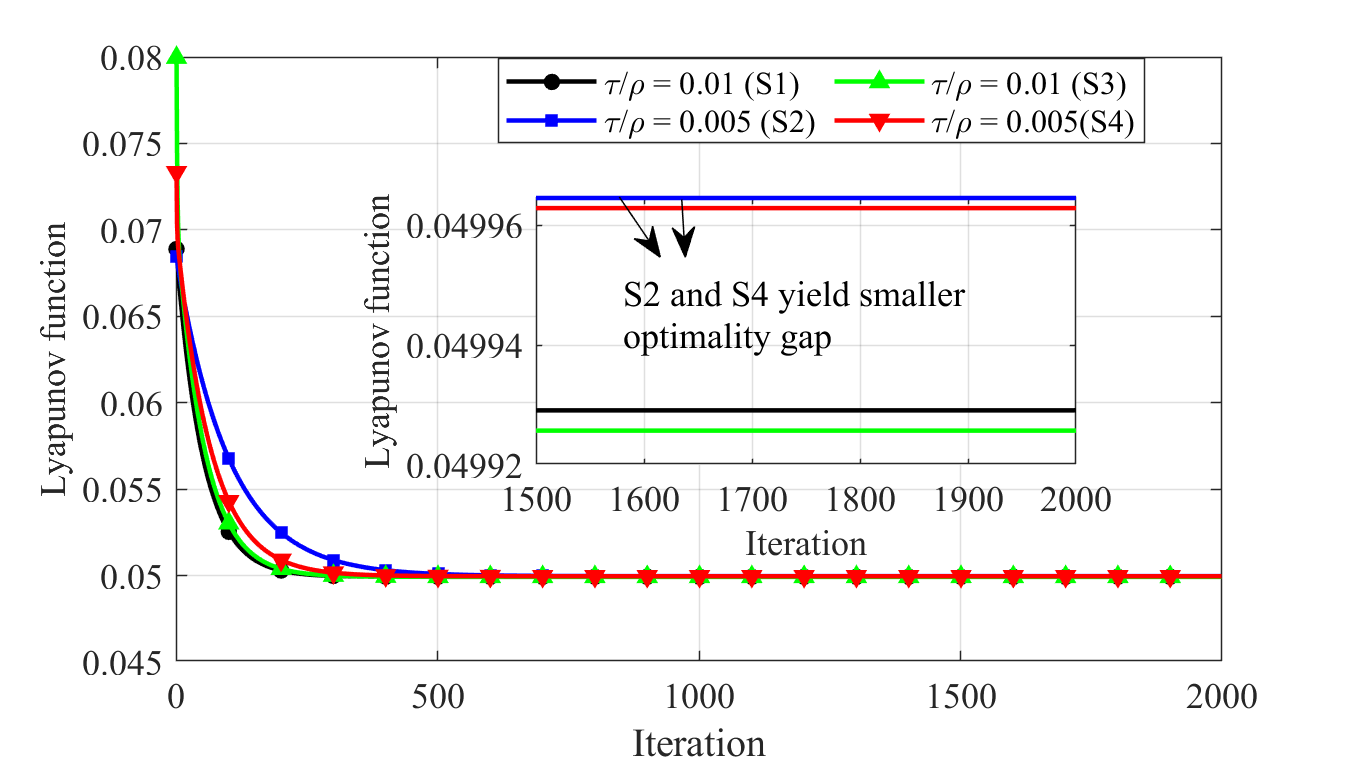}}
 	 	\centerline{\includegraphics[width=3.0 in]{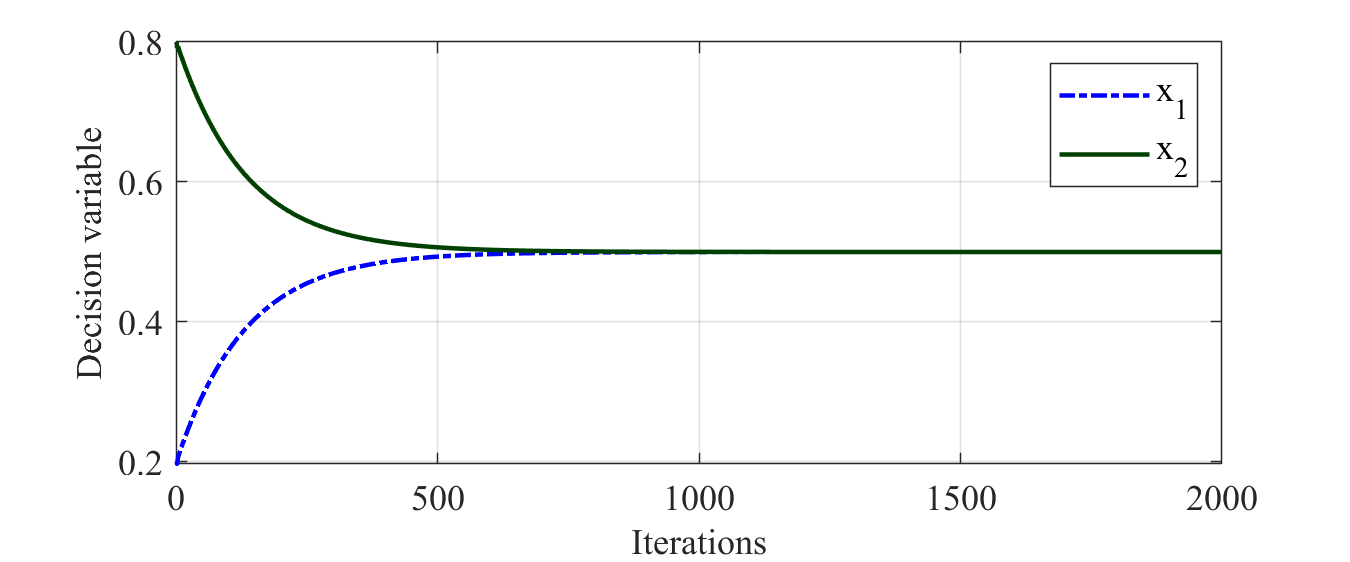}}
 	\caption{(a) The evolution of  Lyapunov function $T_c^k$ with S1-S4. (b) The evolution of primal variables $x_1$ and $x_2$ with S4. }
 	\label{fig:NumericExample}
 \end{figure}

 Before running the algorithm, we first can  easily verify  the convergence condition (C1)  stated in \textbf{Theorem} \ref{thm:theorem1} for S1-S4.  
We use the \texttt{interior-point} method embedded in the \texttt{fmincon} solver of MATLAB to solve  subproblems \eqref{eq:x-update}. 
We run \textbf{Algorithm} \ref{alg: perturbed_ADMM} sufficiently long (i.e., $K=2000$ iterations when the Lyapunov function does not change apparently) for the  settings S1-S4. 
We first examine the convergence of the method indicated by the Lyapunov function.   Fig. \ref{fig:NumericExample} (a) shows the evolution of the Lyapunov function w.r.t. the iterations with S1-S4. 
We observe that for all the settings S1-S4, the Lyapunov functions strictly decrease w.r.t.  the iterations and finally stabilize at some value  that is close to the optima  $f^{\star} = 0.1 x_1^{\star} + 0.1 x_2^{\star} + 0.1 x_1^{\star} x_2^{\star} = 0.05$. By further examining the results, we note that  a larger ratio  $\tau/\rho$  yields faster convergence rate as with S1/S3 ($\tau/\rho= 0.01$) compared with S2/S4 ($\tau/\rho =0.005$).  This is caused by the relatively smaller penalty factor $\rho$ and proximal  factor $\beta$  required to  ensure the convergence condition (C1) for a larger $\tau/\rho$.   Note that the penalty factor $\rho$ and the proximal factor $\beta$ can be interpreted as some means to slow  down the primal updates as they have an effect in penalizing the deviation from current update $\x^k$. 
Oppositely,  a smaller ratio $\tau/\rho$ generally yields higher solution quality (i.e., smaller suboptimality gap) as with S2/S4 ($\tau/\rho = 0.005$) compared with S1/S3 ($\tau/\rho = 0.01$). 
This is in line with \textbf{Theorem} \ref{thm:theorem1}. 

To further examine the solution quality, we report the detailed results with the four settings S1-S4 (\texttt{Prox-ADMM-Sx}, x = 1, 2, 3, 4) and the centralized method   (using the \texttt{interior-point} method embedded in the \texttt{fmincon} solver of MATLAB) in Table \ref{tab:NumericExample}.  Note that the convergent solution $\hat{x}_1$ and $\hat{x}_2$ with proximal ADMM under  the settings S1-S4 are quite close to the optimal solution $x^{\star}_1 = 0.5$ and $x_2^{\star} = 0.5$ obtained from the centralized method. More specifically, by measuring the sub-optimality by   $\Vert \hat{\x} - \x^{\star} \Vert/\Vert \x^{\star} \Vert$ where $\hat{\x} = (\hat{x}_1, \hat{x}_2)$ and $\x^{\star} = (x_1^{\star}, x_2^{\star})$ are the convergent and optimal solution, we have the sub-optimality of  proximal ADMM is around 1.2E-3 with S1/S3 ($\tau/\rho = 0.01$) and 5.9E-4 with S2/S4 ($\tau/\rho = 0.005$).  We therefore imply that a smaller ratio $\tau/\rho$  can achieve higher solution quality  but generally at the cost of slower convergence rate as observed in Fig. \ref{fig:NumericExample} (a). This implies that a trade-off  in terms of the solution quality and the convergence speed is necessary while configuring the algorithm (i.e., the ratio of $\tau/\rho$) for specific applications.
For this example, considering both the solution quality  and convergence rate, we have S4  a preferred option. We therefore display the convergence of  primal variables $x_1$ and $x_2$ with S4 in Fig. \ref{fig:NumericExample}(b). Note that $x_1$ and $x_2$ gradually approach the optimal solution $\x_1^{\star} = 0.5$ and $x_2^{\star} = 0.5$.

\begin{table}[h] 
	\setlength\tabcolsep{1.5pt}
	\centering
	\caption{Performance of proximal ADMM under the settings S1-S4 vs. Centralized method}
	\label{tab:NumericExample}
	\begin{tabular}{lccccc}     
		\toprule
		~~~~Method~~~~ &  $\tau/\rho$ &~~$\hat{x}_1$ &  $\hat{x}_2$ & \makecell{\scriptsize Sub-\\optimality}   &  \makecell{Convergence \\Rate }   \\
	\hline 
		{\scriptsize \texttt{Centralized}} & {--} & {\scriptsize 0.5} &  {\scriptsize 0.5}  &  {--}& {--} \\
		{\scriptsize \texttt{Prox-ADMM-S1}~}  &0.01~ & {\scriptsize 0.4994}  &  {\scriptsize 0.4994} & {\scriptsize1.1E-3}  & No. 1 \\
				{\scriptsize \texttt{Prox-ADMM-S2}~} & 0.005~  & {\scriptsize 0.4997}  &  {\scriptsize 0.4997} & {\scriptsize 5.7E-4} & No. 4 \\
						{\scriptsize \texttt{Prox-ADMM-S3}~} & 0.01~  & {\scriptsize 0.4994}  &  {\scriptsize 0.4994} & {\scriptsize 1.2E-3} & No .2\\
								{\scriptsize \texttt{Prox-ADMM-S4}~} & 0.005~   & {\scriptsize 0.4997}  &  {\scriptsize 0.4997} & {\scriptsize 5.9E-4}  & No. 3\\
		\bottomrule 
	\end{tabular}
\end{table}

\subsection{Application: multi-zone HVAC control}
To showcase the performance of  proximal ADMM in applications, we apply it to the multi-zone heating, ventilation, and air conditioning (HVAC) control arising from smart buildings \cite{yang2021distributed, yang2020hvac, yang2021stochastic}. 
The  goal is to optimize the HVAC operation to provide the comfortable temperature with minimal electricity bill.  Due to the thermal capacity of buildings, the evolution of indoor temperature is a slow process affected both by the dynamic indoor occupancy (thermal loads) and   the HVAC operation (cooling loads).  The general solution is to design a model predictive controller  for optimizing HVAC operation (i.e., zone mass flow and  zone temperature trajectories) to minimize the overall electricity cost while respecting the comfortable temperature ranges based on the predicted information (i.e., indoor occupancy, outdoor temperature, electricity price, etc.). The general problem formulation is presented below. 

\small 
\begin{subequations}
	 		\setlength{\abovedisplayskip}{-2pt}
	\setlength{\belowdisplayskip}{1pt}
		\begin{alignat}{4}
		 \label{pp:HVAC problem} \min_{\m^z, \T} &\sum_{t}  c_t \big\{  c_p  (1-d_r)\sum_{i}  m^{zi}_t (T^o_t-T^c) \tag{$\mathbf{P2}$} \\
		&  + c_p \eta d_r \sum_{i}    m^{zi}_t (T^i_t-T^c) +  \kappa_f  \big(\textstyle\sum_{i}  {m^{zi}_t}\big)^2
		\big\} \Delta_t  \notag \\
	    \text{s.t.}~ &  ~ T^i_{t+1} = A_{ii} T^i_t  + {\textstyle \sum}_{j\in N_i} A_{ij} T^j_t  \notag\\
	& \quad \quad \quad  \label{eq:zone_thermal_dynamics}  + C_{ii}  m^{zi}_t   (T^i_t  -T^c ) +  D^{i}_t, ~\forall i, t. \\
	& \label{eq:comfort} {T}^i_{\min}  \leq  T^i_t  \leq T^i_{\max}, ~~~~\forall i, t. \\ 
 	& \label{eq:zone_mass_flow_rate_limit} m^{zi}_{\min} \leq  m^{zi}_t \leq m^{zi}_{\max}, ~\forall i, t.  \\
 	& \label{eq:total_mass_flow_rate} \textstyle\sum_{i} m^{zi}_t \leq \overline{m}, ~\forall t.  
	\end{alignat}
\end{subequations}
\normalsize
where $i$ and $t$ are zone and time indices, $\T =(T^i_t)_{\forall i, t}$ and $\m^{z} = (m^{zi}_t)_{\forall i, t} $ are zone temperature and the supplied zone mass flow rates, which are decision variables.  Note that we have the HVAC system serves $N$ zones  with thermal couplings (i.e., heat transfer) within a building. 
The other notations are parameters. 
 For example,  $[T^i_{\min},  T^i_{\max}]$ represent the comfortable temperature range of zone $i$.  The problem is subject to the  constraints including  zone thermal couplings \eqref{eq:zone_thermal_dynamics},  comfortable temperature margins \eqref{eq:comfort},  zone mass flow rate limits \eqref{eq:zone_mass_flow_rate_limit}, and  total zone mass flow rate limits \eqref{eq:total_mass_flow_rate}. 
 
 Note that problem \eqref{pp:HVAC problem} is generally in large scale for a commercial building due to the large number of zones and rooms.  This problem represents one of the major challenging problems with smart buildings.  Centralized strategies are generally not suitable due to the computation and communication overheads. In this part, we show how the proximal ADMM can be applied to solve problem \eqref{pp:HVAC problem}  in a distributed manner and thus overcome the computation burden. We first reformulate problem \eqref{pp:HVAC problem} as \eqref{pp:HVAC problem_reformulated}. We have  
 {\small 
 	\begin{subequations}
 		\begin{alignat}{4}
 			& \label{pp:HVAC problem_reformulated} \min_{\m^z, \T} \sum_{t}  c_t \big\{  c_p  (1-d_r)\sum_{i }  m^{zi}_t (T^o_t-T^c) \tag{$\mathbf{P3}$} \\
 			&  + c_p \eta d_r \sum_{i} m^{zi}_t (T^{ii}_t-T^c) +   \kappa_f  \big(\textstyle\sum_{i}  {m^{zi}_t}\big)^2
 			\big\} \Delta_t  \notag \\
 			& + M \textstyle \sum_i  \textstyle \sum_t \big(T^{ii}_{t+1} - A_{ii} T^{ii}_t -   {\sum}_{j\in N_i} A_{ij} T^{ij}_t \notag  \\
 			& \quad \quad - C_{ii}  m^{zi}_t   (T^{ii}_t  -T^c ) - D^{i}_t \big)^2 \notag \\
 			&	\text{s.t.}~ ~  \label{eq:temperature_consensus2} T^{ij}_t = \overline{T}^{j}_t, ~~\forall i, j, t.   \\
 			& \label{eq:comfort2} \quad \quad  T^i_{\min}  \leq  T^{ii}_t  \leq T^i_{\max}, ~~~\forall i, t. \\ 
 			&  \label{eq:comfort_consensus2}  \quad \quad  T^i_{\min}  \leq \overline{T}^i_t \leq  T^i_{\max}, ~~~\forall i, t. \\ 
 			&  \label{eq:zone_mass_flow_rate_limit2} \quad \quad m^{zi}_{\min} \leq  m^{zi}_t \leq m^{zi}_{\max}, ~~~\forall i, t.  \\
 			& \label{eq:total_mass_flow_rate2} \quad  \quad \textstyle \sum_{i} m^{zi}_t \leq \overline{m}, ~\forall t.  
 		\end{alignat}
 \end{subequations}}
 where we have augmented the decision variables for each zone controller to involve the copy of temperature for its neighboring zones, i.e., $\T^i: = (T^{ij}_t)_{j \in N_i, ~t}$ with $T^{ij}_t$ denoting the estimated temperature of zone $j$ by zone $i$ and $N_i$ representing the neighboring zones of zone $i$. Besides, to drive the consistence of duplicated zone temperature copies, we introduce a block of consensus variable $\overline{\T} = (\overline{T}^i_t)_{i, t}$. Considering the challenging to handle the hard non-linear constraints \eqref{eq:zone_thermal_dynamics}, we turn to the penalty method and penalize the violations by quadratic terms. 
Note that problem \eqref{pp:HVAC problem_reformulated} fits into the template \eqref{pp:problem template} and we have $N + 1$ computing agents, where agents $1$ to $N$ correspond to the zones with the augmented  decision variable $\x_i = ( (T^{ij}_t)_{j \in N_i, t}, (m^{zi}_t)_{t})$, and agent $0$ control the consensus decision variable $\overline{\mathbf{T}}= (\overline{T}_i)_{i \in N}$. Constraints \eqref{eq:temperature_consensus2} and \eqref{eq:total_mass_flow_rate2} represents the coupled linear constraints which can be expressed in the compact form $\A\x = \bb$ if necessary. The other constraints comprise the local bounded convex constraints $\X_i$ for the agents.

\begin{figure}[h]
	\centerline{\includegraphics[width=3.1 in]{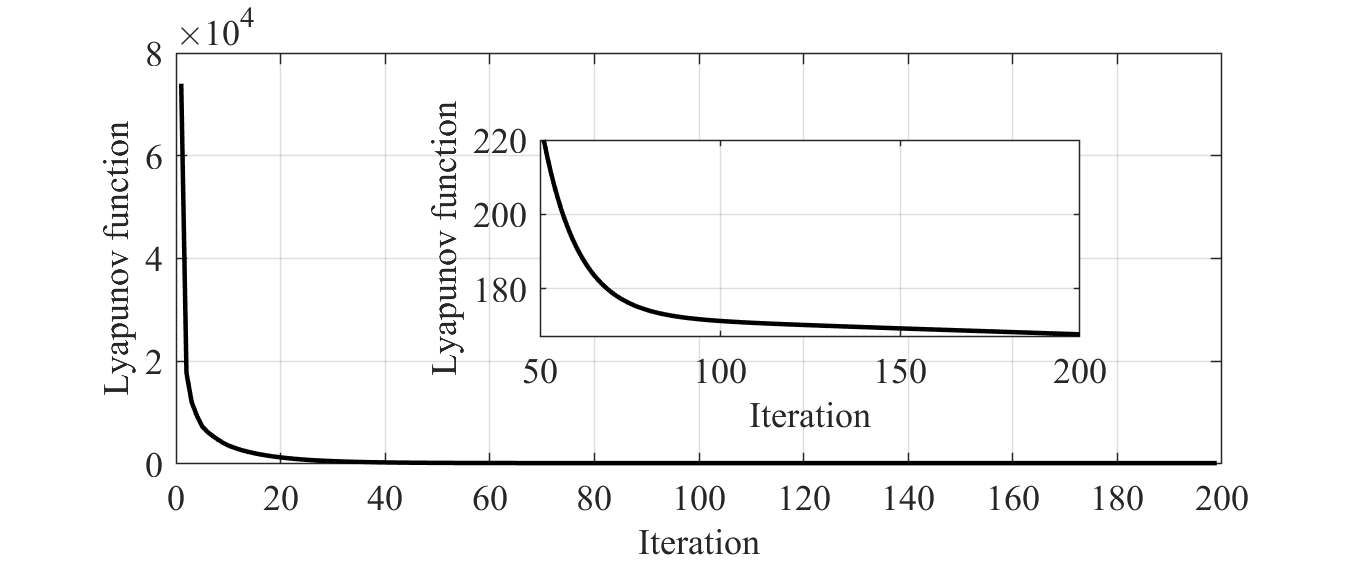}}
		\centerline{\includegraphics[width=3.1in]{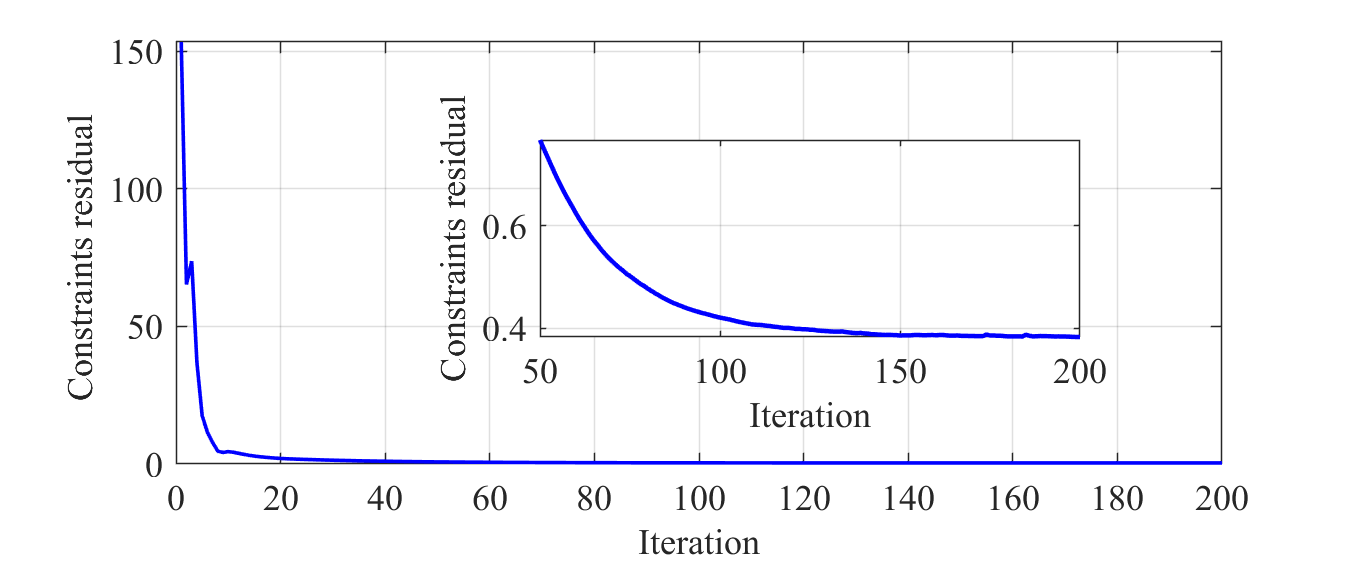}}
	\caption{(a) The evolution of Lyapunov function. (b) The evolution of the norm of constraints residual. }
	\label{Fig:HVAC_potential_residual}
\end{figure}

We consider a case study with $N = 10$ zones and the predicted horizon $T = 48$ time slots (a whole day with a sampling interval of 30 \si{mins}).  We set the lower and upper comfortable temperature bounds as $T^i_{\min} = 24^{\circ}$C and $T^i_{\max} = 26^{\circ}$C.  The specifications for  HVAC system can refer to \cite{yang2020hvac, yang2021distributed}. 
The algorithm of  proximal ADMM is configured by  $\rho = 2.0$, $\tau = 0.1$, $\beta = 3.0$,  $\B_i = \mathbf{I}$ (suitable sizes), and $c = 8.7$.   We run the algorithm suitably long  ($K =200$ iterations when  both the residual and Lyapunov function do not change apparently). We first examine the  convergence of the algorithm measured by the \texttt{Lyapunov  function} and the norm of (coupled) \texttt{constraints residual}.
We visualize   the \texttt{Lyapunov function} and  \texttt{constraints residual}  in Fig. \ref{Fig:HVAC_potential_residual}. 
Note that the \texttt{Lyapunov function}  strictly declines along the iterations, which is consistent with our theoretical analysis.  Besides, the \texttt{constraints residual}  almost  strictly decreases  with the iterations as well and finally approaches \emph{zero}. 
We have the overall norm of the \texttt{constraints residual} at the end of  iterations is about $0.38$, which is quite small considering the problem scale $T \cdot N = 480$.  This justifies  the convergence 
property of  proximal ADMM for the smart building application.

\begin{figure}[h]
  		\centerline{\includegraphics[width=3.0 in]{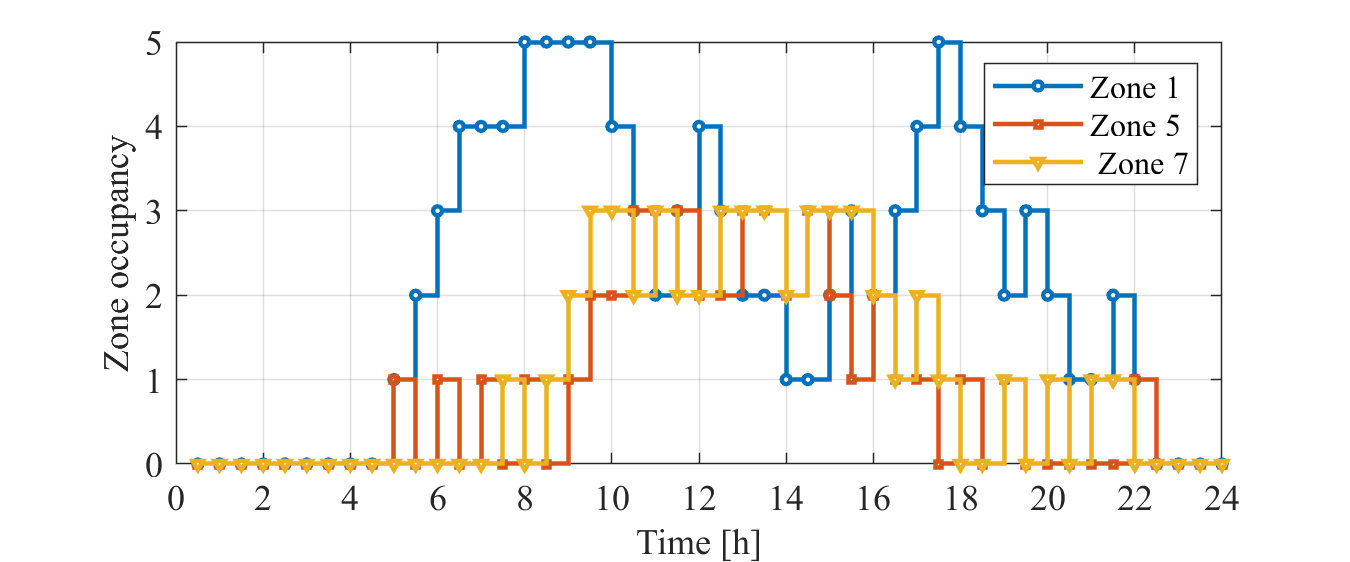}}
 	\centerline{\includegraphics[width=3.0 in]{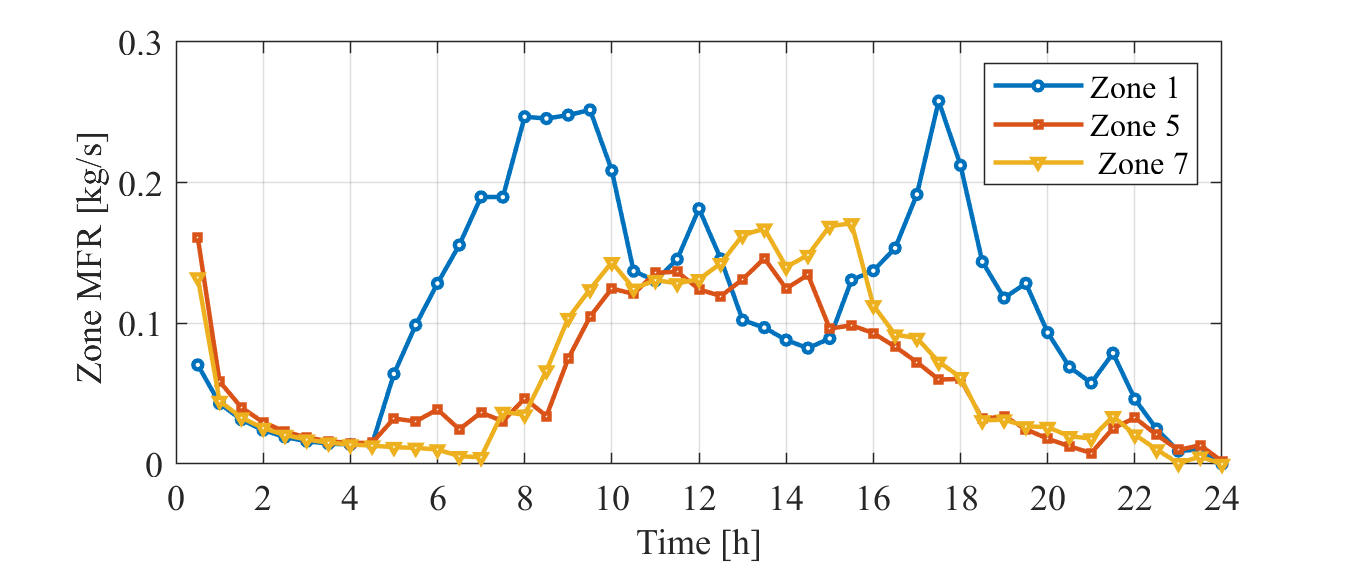}}
 		\centerline{\includegraphics[width=3.0 in]{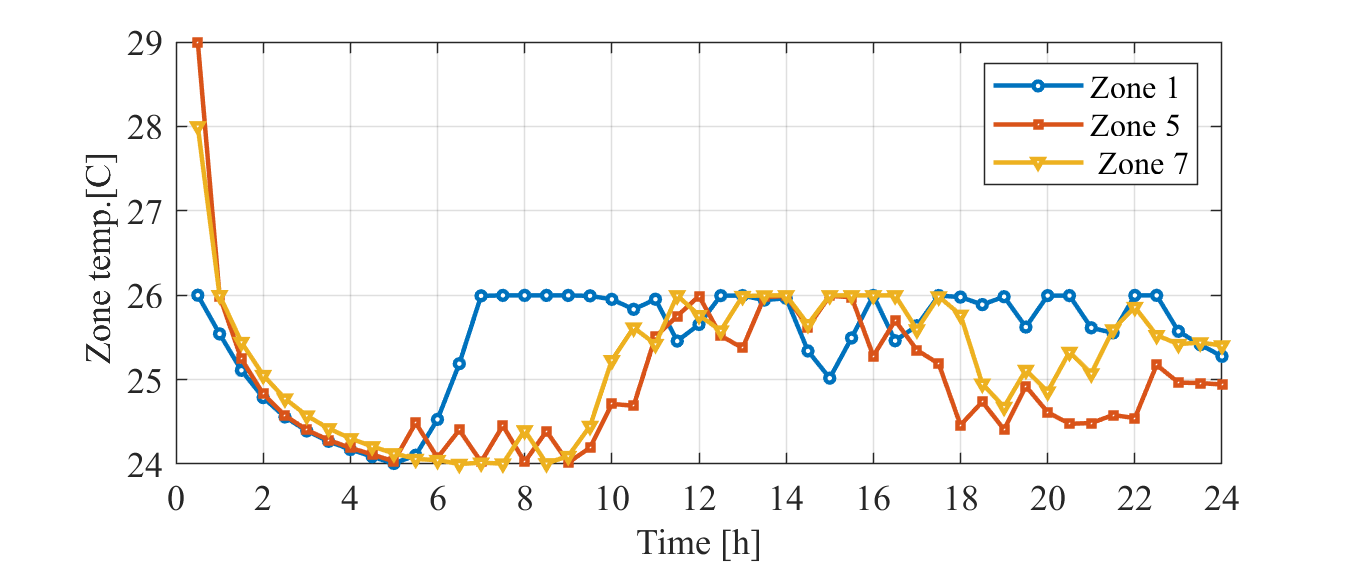}}
 	\caption{(a) Zone occupancy. (b) Zone mass flow rate (Zone MFR). (c) Zone temperature (Zone temp.). }
 	\label{Fig:HVAC_occupancy_temperature_MFR}
 \end{figure}

We next evaluate the solution quality measured by the HVAC electricity cost and human comfort.  We randomly pick $3$ zones (zone 1, 3, 7) and display the predicted \texttt{zone occupancy} (inputs), the zone mass flow rates (\texttt{zone MFR}, control variables) ,  and the zone temperature (\texttt{zone temp.}, control variables) over the $48$  time slots in Fig. \ref{Fig:HVAC_occupancy_temperature_MFR}.  Note that the variations of  \texttt{zone MFR}  are  almost consistent with the \texttt{zone occupancy}. This is reasonable as the \texttt{zone occupancy}  determines the thermal loads which need to be balanced by the zone mass flow rates.  
We besides see that the \texttt{zone temp.} are all maintained within the comfortable range $[24, 26]^{\circ}$C.  This infers  the satisfaction of  \texttt{human comfort}. 
To further evaluate the solution quality and computation efficiency, we compare the proximal ADMM (\texttt{Prox-ADMM}) with centralized method (\texttt{Centralized}). Specifically, we use the \texttt{interior-point} embedded in the \texttt{fmincon} solver of MATLAB to solve both the subproblems \eqref{eq:x-update} with \texttt{Prox-ADMM} and   problem \eqref{pp:HVAC problem} with \texttt{Centralized}. 
For the \texttt{Centralized}, we run the solver sufficiently long without considering the  time with the objective to approach the best possible optimal solution.
 We compare the two methods  in three folds, i.e.,  \texttt{electricity cost}, the norm of \texttt{constraints residual},  and \texttt{computation time} as reported in Table \ref{tab:HVAC}. 
We see that  \texttt{electricity cost}  with  \texttt{Prox-ADMM} is about $160.20$ (s\$) versus $153.12$ (s\$) yield by \texttt{Centralized}.  We imply the sub-optimality of  \texttt{Prox-ADMM} in terms of the objective  is about  $5.0\%$. 
However, the \texttt{Prox-ADMM}
obviously outperforms the \texttt{Centralized} in  computation efficiency.  The average \texttt{computing time} for each zone is about 50~\si{min} with \texttt{Prox-ADMM} (parallel computation) while the \texttt{Centralized} takes more than 10~\si{\hour}. Note that we have picked $T=48$ time slots (a whole day) as the predicted horizon, the computing time could be largely sharpened in practice with a much smaller prediction horizon, say $T=10$ time slots (5\si{\hour}).

 \begin{table}[h] 
 	\setlength\tabcolsep{2pt}
 	\centering
 	\caption{Prox-ADMM vs. Centralized  for HVAC control in smart buildings ($N=10$ zones)}
 	\label{tab:HVAC}
 	\begin{tabular}{lcccc}     
 		\toprule
 		~~Method~ & \makecell{ Electricity \\ cost (s\$)} & \makecell{Human \\ comfort} &  \makecell{Constraints \\residual}  &   \makecell{Computing \\ time }\\
 		\hline 
 		\texttt{Centralized}  & 153.12 &  Y & 0 & $\geq 10$\si{\hour}\\
 		\texttt{Prox-ADMM}    & 160.54 &  Y & 0.38  & 50~\si{\minute}  \\
 		\bottomrule 
 	\end{tabular}
 \end{table}

 \section{Conclusion and Future Work}
This paper focused on developing a distributed algorithm for a class of  nonconvex and nonsmooth problems with convergence guarantee.  
The problems  are featured by  i) a possibly nonconvex objective  composed of both separate and composite components,  ii) local bounded convex constraints,  and iii) global coupled linear constraints.  This class of problems is broad in application but lacks distributed methods with convergence guarantee. 
We turned to the powerful alternating direction method of multiplier (ADMM)  for constrained optimization 
but faced the challenge to establish convergence. 
Noting that the underlying obstacle is to assume the boundness of dual updates, we revised the classic ADMM and proposed to update the dual variables in a distributed manner. This leads to a proximal ADMM with  the convergence guarantee towards the approximate stationary points of the problem.  We demonstrated  the convergence and solution quality of  the distributed method by  a numerical example and  a concrete  application to the multi-zone heating, ventilation, and air-condition (HVAC) control arising from smart buildings. 

 This paper proposed the discounted dual update scheme in conjunction with ADMM for a class of nonconvex and nonsmooth problems, some interesting future work includes studying whether the discounted dual update scheme can be explored to develop distributed methods for more broad classes of problems both in convex and  nonconvex settings.
 

\bibliographystyle{ieeetr}        
\bibliography{reference}          

\appendix
\section{Proof of Proposition \ref{prop:prop-x}}\textbf{Prop.} \ref{prop:prop-x} is  established based on the first-order optimality condition of subproblems \eqref{eq:x-update} and the $\Lips$ continuous gradient property of  $f$ and $g$. 

We first establish the following equality and notation. 
{\small 
			 		\setlength{\abovedisplayskip}{-2pt}
	\setlength{\belowdisplayskip}{1pt}
\begin{align} 
	\label{eq:equality}&\A_i\x^{k+1}_i  + {\textstyle \sum}_{j \neq i} \A_j \x^k_j - \bb  \\
	& ~=\! \A \x^k -\bb + \A_i(\x_i^{k+1}-\x_i^k). \notag\\
	&~= \!\A \x^{k+1} \!-\!\bb + \A (\x^k - \x^{k+1}) + \A_i(\x_i^{k+1}-\x_i^k). \notag\\
	\label{eq:hatlambda}&\hat{\boldlambda}^k:  = \boldlambda^{k} + \rho (\A \x^{k + 1} - \bb).
\end{align}
}


For  subproblems \eqref{eq:x-update}, the first-order optimality condition states  that there exists $\nu_{i}^{k + 1} \in N_{\X_i}(\x_{i}^{k + 1})$ that
{\small 	
	\setlength{\abovedisplayskip}{-2pt}
	\begin{align*} %
		0  & =   ~\nabla f_i(\x_i^{k+1}) + \nabla_i g(\x^{k}) +  \A^\top_i \boldlambda^k \\
		& \quad \quad + \rho \A^\top_i ( \A_i\x^{k+1}_i  + {\textstyle \sum}_{j \neq i} \A_j \x^k_j  - \bb)  \notag\\ 
		&\quad \quad  + \beta \B^\top_i \B_i( \x^{k+1}_i-\x^k_i) + \nu_{i}^{k + 1}   \notag\\
		&  = ~\nabla  f_i(\x_i^{k+1})  +  \nabla_{i} g(\x^{k}) + \A^\top_i  \big( \boldlambda^k +  \rho  ( \A \x^{k+1}  - \bb ) \big) \notag\\
		& \quad  \quad +  \rho\A_i^\top  \A (\x^k - \x^{k+1}) +  \rho \A_{i}^{\top} \A_i(\x_i^{k+1} - \x_i^k)  \\
		& \quad  \quad +  \beta \B^\top_i \B_i ( \x^{k+1}_i-\x^k_i) + \nu_{i}^{k + 1} ~~~~~~~~~~\text{by } \eqref{eq:equality} \notag\\
		&  = ~\nabla  f_i(\x_i^{k+1}) + \nabla_{i} g(\x^{k})  +\A^\top_i \hat{\boldlambda}^{k} + \rho \A^\top_i    \A (\x^k - \x^{k+1}) \notag\\
		&  \quad  \quad + \rho \A^\top_i \A_i(\x_i^{k+1}-\x_i^k) \notag\\
		& \quad  \quad  + \beta  \B^\top_i\B_i\big( \x^{k+1}_i-\x^k_i\big)  + \nu_{i}^{k + 1}~~~~~~~~~~\text{by}~\eqref{eq:hatlambda}. \notag
	\end{align*}
}
%
Multiplying by ($\x_{i}^{k + 1} - \x_{i}$) in both sides, we have
{\small 
\begin{equation} \label{eq:first-order-equations}
		\setlength{\abovedisplayskip}{-2pt}
	\setlength{\belowdisplayskip}{1pt}
	\begin{split}
		&\langle \nabla    f_i(\x_i^{k+1}),   \x_{i}^{k + 1} - \x_{i}\rangle + \langle \nabla_{i} g(\x^{k}), \x_{i}^{k + 1} - \x_{i}\rangle \\
		& + \langle \hat{\boldlambda}^{k}, \A_i(\x_{i}^{k + 1} - \x_{i}) \rangle \\
		&+ \rho \langle \A (\x^k - \x^{k+1}),  \A_{i}(\x_{i}^{k + 1} - \x_{i})\rangle\\\
		&  + \rho \langle \A_i(\x_i^{k+1}-\x_i^k),  \A_{i}(\x_{i}^{k + 1} - \x_{i})\rangle\\
		&  + \beta  \langle \B_i( \x^{k+1}_i-\x^k_i), \B_i(\x_{i}^{k + 1} - \x_{i})\rangle \\
		= &   -\langle \nu_{i}^{k + 1},  \x_{i}^{k + 1} - \x_{i} \rangle \leq 0, ~~\forall \x_{i} \in \X_i. 
	\end{split}
\end{equation}}

Summing  up \eqref{eq:first-order-equations} over $i$,  we have $\forall \x_{i} \in \X_{i}$,
{\small 
	\setlength{\abovedisplayskip}{-2pt}
	\begin{align*} 
		\langle \nabla  & f(\x^{k+1}), \x^{k+1}- \x\rangle  
		+ \langle \nabla  g(\x^{k}), \x^{k+1}- \x\rangle \\
		& + \langle  \hat{\boldlambda}^{k}, \A(\x^{k+1}-\x)\rangle + (\x^{k+1}-\x)\rho\A^\top\A (\x^k-\x^{k+1})\\
		& + \textstyle \sum_{i}(\x^{k+1}_i-\x_i)^\top (\rho\A^\top_i\A_i + \beta\B^\top_i\B_i)(\x^{k+1}_i-\x^k_i) \leq 0.  \\
\end{align*} }

Plugging in $\Q:= \rho G_{\A} + \beta G_{\B} - \rho\A^\top\A$, we  have 
{\small 
	\setlength{\abovedisplayskip}{-2pt}
	\begin{align} \label{eq:first-order-k}
		\langle \nabla  & f(\x^{k+1}), \x^{k+1}- \x\rangle  
		+ \langle \nabla  g(\x^{k}), \x^{k+1}- \x\rangle \notag\\
		& + \langle  \hat{\boldlambda}^{k}, \A(\x^{k+1}-\x)\rangle  \notag\\
		& + (\x^{k+1}-\x)^\top \Q (\x^{k+1}-\x^k) \leq 0,  ~~\forall \x \in \X.  
	\end{align}
}

By induction, we have 
{\small 
	\setlength{\abovedisplayskip}{-2pt}
	\setlength{\belowdisplayskip}{1pt}
\begin{align} \label{eq:first-order-k-1}
	 \langle \nabla &f(\x^k), \x^k- \x \rangle
	+ \langle \nabla  g(\x^{k-1}), \x^{k}- \x\rangle \notag\\
	& + \langle \hat{\boldlambda}^{k-1}, \A(\x^k-\x)\rangle  \notag\\
	& + (\x^k-\x)^\top \Q (\x^k-\x^{k-1}) \leq 0, ~~\forall \x \in \X.
\end{align}}

By setting $\x = \x^k$ and $\x=\x^{k+1}$  with \eqref{eq:first-order-k} and  \eqref{eq:first-order-k-1},  we have
{\small 
	\setlength{\abovedisplayskip}{-2pt}
	\begin{align}
		\langle \nabla & f(\x^{k+1}), \x^{k+1}- \x^{k}\rangle  
		+ \langle \nabla  g(\x^{k}), \x^{k+1}- \x^{k}\rangle \notag\\
		& + \langle  \hat{\boldlambda}^{k}, \A(\x^{k+1}-\x^{k})\rangle  \notag\\
		&  \label{eq:k}+ (\x^{k+1}-\x^{k})^\top \Q (\x^{k+1}-\x^k) \leq 0.\\
		\langle \nabla &f(\x^k), \x^k- \x^{k+1} \rangle
		+ \langle \nabla  g(\x^{k-1}), \x^{k}- \x^{k + 1}\rangle \notag\\
		& + \langle \hat{\boldlambda}^{k-1}, \A(\x^k-\x^{k + 1})\rangle  \notag\\
		&\label{eq:k-1} + (\x^k-\x^{k + 1})^\top \Q (\x^k-\x^{k-1}) \leq 0.
	\end{align}
}

Summing up \eqref{eq:k} and \eqref{eq:k-1} and  plugging  in  $\w^k: = (\x^{k+1}-\x^k) - (\x^k-\x^{k-1})$, we have 
\begin{equation} \label{eq:first-order-sum}
	\begin{split}
		\langle \nabla & f(\x^{k+1}) - \nabla  f(\x^k), \x^{k+1}- \x^k \rangle \\
		& +  \langle \nabla g(\x^{k}) - \nabla g(\x^{k-1}), \x^{k+1}- \x^k \rangle \\
		&  + \langle \hat{\boldlambda}^{k} - \hat{\boldlambda}^{k-1}, \A(\x^{k+1} -\x^k)\rangle \\
		& + (\x^{k+1}-\x^k)^{\top}\Q\w^k  \leq  0. 
	\end{split}
\end{equation}  

Based on the $\Lips$ continuous gradient property of $f$ over the compact set $\x \in \X$, we have 

{\small
	\begin{equation} \label{eq:first-term}
		\begin{split}
			\langle  \nabla f(\x^{k+1})\! -\! \nabla f(\x^k), \x^{k+1} \!-\! \x^k\rangle\! \geq \!- L_f \Vert \x^{k + 1} \!-\! \x^k  \Vert^2. 
		\end{split}
	\end{equation}
}

We also have
{\small 
	\begin{align*}
		& \langle \nabla  g(\x^{k}) - \nabla g(\x^{k-1}), \x^{k + 1} - \x^{k}\rangle  \\
		& =\langle \frac{\nabla g(\x^{k}) - \nabla g(\x^{k-1})}{\sqrt{L_{g}}}, \sqrt{L_{g}}(\x^{k + 1} - \x^{k})\rangle \\
		& \geq -\frac{1}{2L_{g}} \Vert \nabla g(\x^{k}) - \nabla g(\x^{k-1})\Vert^{2} - \frac{L_{g}}{2}\Vert \x^{k + 1} - \x^{k}\Vert^{2} \\
		& \geq -\frac{L_{g}}{2} \Vert \x^{k} - \x^{k-1} \Vert^{2} - \frac{L_g}{2} \Vert \x^{k + 1} - \x^{k} \Vert^{2}
\end{align*}}
where the last equality is based on  the $\Lips$ continuous gradient property of $g$. 

Besides, we have 
{\small
\begin{equation} \label{eq:second-term}
	\begin{split}
		& \quad \langle \hat{\boldlambda}^{k}- \hat{\boldlambda}^{k-1}, \A(\x^{k + 1} - \x^{k})\rangle \\
		& = \left\langle \boldlambda^{k+1} - \boldlambda^k + \tau(\boldlambda^k - \boldlambda^{k-1}), \A(\x^{k+1} - \x^k)\right \rangle \\
		& = \Big\langle  \boldlambda^{k+1} - \boldlambda^k + \tau(\boldlambda^k - \boldlambda^{k-1}),   \\
		& \quad \quad \quad \quad \frac{\boldlambda^{k+1}   - \boldlambda^k}{\rho}  - \frac{(1-\tau)}{\rho}(\boldlambda^k-\boldlambda^{k-1})\Big \rangle \\
		& = \frac{\Vert \boldlambda^{k+1} \!-\! \boldlambda^k \Vert^2}{\rho} \!-\! \frac{(1-2\tau)}{\rho} \langle \boldlambda^{k+1} - \boldlambda^k, \boldlambda^k - \boldlambda^{k-1} \rangle \\
		& \quad  -\frac{ \tau(1-\tau)}{\rho}\Vert \boldlambda^k - \boldlambda^{k-1} \Vert^2  \\
		& \geq  \frac{\Vert \boldlambda^{k+1} - \boldlambda^k \Vert^2}{\rho}  -\frac{1-2\tau}{2\rho} \Vert \boldlambda^{k+1} - \boldlambda^k \Vert^2  \\
		&  \quad  -  \frac{1-2\tau}{2\rho}\Vert \boldlambda^k - \boldlambda^{k-1}\Vert^2  - \frac{\tau(1-\tau)}{\rho}\Vert \boldlambda^k - \boldlambda^{k-1} \Vert^2 \\
		& =  \frac{1-2\tau^2}{2\rho} \Vert \boldlambda^{k+1} - \boldlambda^k \Vert^2  -  \frac{1-2\tau^2}{2\rho} \Vert \boldlambda^k - \boldlambda^{k-1} \Vert^2 \\
		& \quad  +{ \tau(\tau + 1)}/{\rho} \Vert\boldlambda^{k+1} - \boldlambda^k  \Vert^2
	\end{split}
\end{equation}}
where the inequality is by $\langle \mathbf{a}, \mathbf{b} \rangle \leq \frac{1}{2}(\Vert \mathbf{a} \Vert^{2} + \Vert \mathbf{b}\Vert^{2})$.

Based on the inequality $ \mathbf{b}^{\top} \mathbf{M} (\mathbf{b} - \mathbf{c}) = \frac{1}{2}(\Vert \mathbf{b} - \mathbf{c} \Vert_{\mathbf{M}}^{2} + \Vert \mathbf{b} \Vert_{\mathbf{M}}^{2} - \Vert  \mathbf{c} \Vert_{\mathbf{M}}^{2})$, and by setting  $\mathbf{M} = \Q$,  $\mathbf{b}= \x^{k + 1} - \x^{k}$,  and $\mathbf{c} = \x^{k} - \x^{k-1}$,  we have
{\small 
	\begin{equation} \label{eq:third-term}
		\setlength{\abovedisplayskip}{-2pt}
		\begin{split}
			(\x^{k+1}-\x^k)^{\top}\Q\w^k = & \frac{1}{2}  (  \Vert  \w^k \Vert^2_{\Q}  + \Vert  \x^{k+1}-\x^k \Vert^2_{\Q} \\
			&- \Vert \x^k-\x^{k-1} \Vert^2_{\Q} ).
		\end{split}
	\end{equation}
}

Plugging  \eqref{eq:first-term}, \eqref{eq:second-term}, \eqref{eq:third-term} into \eqref{eq:first-order-sum},  we have
{\small
	\begin{align*}
		& \frac{1-2\tau^2}{2\rho}   \left\Vert \boldlambda^{k+1} -\boldlambda^k \right\Vert^2   + \frac{1}{2}\Vert  \x^{k+1}-\x^k\Vert^2_{\Q} \\
		& \quad + \frac{L_g}{2}\Vert \x^{k +1} - \x^{k}\Vert^{2} + \frac{1}{2} \Vert \w^k \Vert^2_{\Q}\notag \\
		&  \leq  \frac{1-2\tau^2}{2\rho}   \left \Vert  \boldlambda^k-\boldlambda^{k-1} \right\Vert^2 + \frac{1}{2}\Vert  \x^k-\x^{k-1}\Vert^2_{\Q} \\
		&\quad + \frac{L_g}{2}\Vert \x^{k} - \x^{k - 1}\Vert^{2} + (L_g + L_f)\left\Vert \x^{k+1} -\x^k \right\Vert^2 \\
		& \quad- {\tau(1 + \tau)}/{\rho} \Vert \boldlambda^{k+1}-\boldlambda^k\Vert^2.
\end{align*}}
We therefore complete the proof. 

\section{Proof of Proposition \ref{prop:prop-L}}
Before starting the proof, we first establish the following inequalities to be used. 
Based on the $\Lips$ continuous gradient property of $f: \R^n \rightarrow \R$ over $\x \in \X$ (see (A1)), we have \cite{guo2017convergence}
\begin{align} \label{eq:Lipsf}
	 f(\x^{k+1}) - f(\x^k) \leq ~ &\langle \nabla f(\x^{k+1}), \x^{k+1}- \x^k \rangle \notag \\
	& + L_f/2 \Vert \x^k - \x^{k+1} \Vert^2. 
\end{align}

Similarly, for $g: \R^n \rightarrow \R$ with $\Lips$ continuous gradient  over $\x \in \X$ (see (A1)),  we have \cite{guo2017convergence}
\begin{align} \label{eq:Lipsg}
	g(\x^{k+1}) - g(\x^k) \leq ~ &\langle \nabla g(\x^{k}), \x^{k+1}- \x^k \rangle \notag \\
	&+ L_g/2 \Vert \x^k - \x^{k+1} \Vert^2. 
\end{align}

Besides, we have
\begin{align} 
	&~~~\frac{\rho}{2}\Vert \A\x^{k+1} - \bb \Vert^2 -  \frac{\rho}{2}\Vert \A\x^k - \bb \Vert^2  \notag\\
	= &~~~\frac{\rho}{2} \left\langle  \A(\x^{k+1} - \x^k),  \A\x^{k+1} + \A\x^k - 2\bb \right\rangle  \label{eq:inequality2} \\
	=&~~~\frac{\rho}{2} \left\langle  \A(\x^{k+1} - \x^k),  \A(\x^k - \x^{k+1}) + 2(\A\x^k - \bb) \right\rangle  \notag\\
	= & -\!\frac{\rho}{2} \Vert \A(\x^{k+1} \!\!-\!\! \x^k)\Vert^2 \!+\! \left\langle \A(\x^{k+1} \!\!-\!\! \x^k), \rho(\A\x^{k+1}\!\!-\!\!\bb)\right\rangle. \notag
\end{align}

We next quantify the decrease of   $\Lag_{\rho}(\x, \boldlambda)$ with respect to (w.r.t.) the primal updates. We have 
{\small
\begin{equation} \label{eq:L-x-decrease}
	\begin{split}
		& \Lag_{\rho}(\x^{k+1}, \boldlambda^k)   - \Lag_{\rho}(\x^k, \boldlambda^k)  \\
		=&~ f(\x^{k+1}) \!-\! f(\x^k) + g(\x^{k + 1}) \!-\! g(\x^{k}) + \langle \boldlambda^k, \A(\x^{k+1}-\x^k) \rangle  \\
		& \quad  + \frac{\rho}{2}\Vert \A\x^{k+1} - \bb \Vert^2 -  \frac{\rho}{2}\Vert \A\x^k - \bb \Vert^2 \\
		\leq  &  ~ \langle \nabla  f(\x^{k+1}) + \nabla g(\x^{k}), \x^{k+1} - \x^k \rangle + \rho_F/2 \Vert \x^k - \x^{k+1} \Vert^2 \\
		& \quad +  \left\langle\boldlambda^k, \A(\x^{k+1}-\x^k)  \right\rangle  -\frac{\rho}{2} \Vert \A(\x^{k+1} - \x^k)\Vert^2 \\
		&  \quad \!+\! \left\langle \A(\x^{k+1} - \x^k), \rho(\A\x^{k+1}-\bb)\right\rangle ~\text{by \eqref{eq:Lipsf}, \eqref{eq:Lipsg}, \eqref{eq:inequality2}} \\
		= & ~ \langle \nabla  f(\x^{k+1}) +  \nabla g(\x^{k}) + \A^\top \hat{\boldlambda}^{k}, \x^{k+1} - \x^k \rangle \\
		& \quad  +  \rho_F/2 \Vert \x^k - \x^{k+1} \Vert^2 - {\rho}/{2} \Vert \A(\x^{k+1} - \x^k)\Vert^2  ~~\text{by \eqref{eq:hatlambda}}\\
		\leq & - \!\Vert \x^{k+1} \!-\! \x^k \Vert_{\Q}^2 \!+\! {\rho_F}/{2} \Vert \x^k \!-\! \x^{k+1} \Vert^2\!\\
		& \quad  -\!{\rho}/{2}~ \Vert \A(\x^{k+1} - \x^k)\Vert^2  ~~\text{by \eqref{eq:k}}.\\ 
	\end{split}	
\end{equation}}

We next  quantify the change of  $\Lag_\rho(\x, \boldlambda)$ w.r.t.  dual update. We have 
\begin{equation} \label{eq:L-lambda-decrease}
	\begin{split}
		& \Lag_\rho(\x^{k+1}, \boldlambda^{k+1}) - \Lag_\rho(\x^{k+1}, \boldlambda^k) \\
		=& \left\langle \boldlambda^{k+1} -\boldlambda^k, \A\x^{k+1}-\bb\right\rangle \\
		= &  \left\langle \boldlambda^{k+1} -\boldlambda^k,  \frac{\boldlambda^{k+1} - (1-\tau)\boldlambda^k}{\rho}\right\rangle  \\
		= & \left\langle \boldlambda^{k+1} \!-\!\boldlambda^k,  \frac{1-\tau}{\rho}(\boldlambda^{k+1}-\boldlambda^k) +\frac{\tau}{\rho} \boldlambda^{k+1}\right\rangle  \\
		= & \frac{(1-\tau)}{\rho}\Vert  \boldlambda^{k+1} - \boldlambda^k\Vert^2 + \frac{\tau}{2\rho}\Big( \Vert \boldlambda^{k+1} - \boldlambda^k \Vert^2 \\
		&  + \Vert \boldlambda^{k+1} \Vert^2 - \Vert \boldlambda^k \Vert^2\Big) \\
		= & \frac{2-\tau}{2\rho} \Vert \boldlambda^{k+1} \!-\! \boldlambda^k\Vert^2 \!+\! \frac{\tau}{2\rho}\Vert \boldlambda^{k+1} \Vert^2 -  \frac{\tau}{2\rho}\Vert \boldlambda^k \Vert^2.
	\end{split}
\end{equation} 


Combining \eqref{eq:L-x-decrease} and \eqref{eq:L-lambda-decrease}, we have 
\begin{equation*} \label{eq:regularized_AL}
	\begin{split}
		& \Lag_{\rho}(\x^{k+1}, \boldlambda^{k + 1})  \!-\!  \frac{\tau}{2\rho}\Vert \boldlambda^{k+1} \Vert^2 -  \big( \Lag_{\rho}(\x^k, \boldlambda^k)  \!-\! \frac{\tau}{2\rho}\Vert \boldlambda^k \Vert^2 \big)\\
		&  \leq  - \Vert \x^{k+1} - \x^k \Vert_{\Q}^2 + \frac{\rho_F}{2} \Vert \x^k -\x^{k+1} \Vert^2\!\\
		& \quad  -\frac{\rho}{2} \Vert \A(\x^{k+1} - \x^k)\Vert^2  + \frac{2-\tau}{2\rho} \Vert \boldlambda^{k+1} - \boldlambda^k\Vert^2.
	\end{split}
\end{equation*}
We have $\Lag_{\rho}^{+}(\x, \boldlambda) = \Lag_{\rho}(\x, \boldlambda) - \frac{\tau}{2\rho} \Vert \boldlambda\Vert^2$, we therefore close the proof. 
\section{Proof of Corollary \ref{cor:corollary}}    
	i) Prove $\Vert \boldlambda^{*} \Vert^2 \leq  \rho \tau^{-1} T_c^0$: Based on the sufficiently decreasing property of $T_c^{k+1}$ (see \textbf{Prop.} \ref{prop:sufficient-decrease}), we have 
		{\small 
		\setlength{\abovedisplayskip}{-4pt}
	\begin{align}\label{eq:sufficient_decrease} 
		T_c^{k+1} \leq T_c^{0}
	\end{align}}

Recalling the definition of the Lyapunov function in \eqref{eq:T-def} and invoking \eqref{eq:second-term_L}, we have
	{\small 
		\setlength{\abovedisplayskip}{-2pt}
		\setlength{\belowdisplayskip}{1pt}
		\begin{equation} \label{eq:T_def2}
			\begin{split}
				& T_c^{k+1} = f(\x^{k+1}) + g(\x^{k + 1})  + \frac{\tau}{\rho} \Vert \boldlambda^{k+1} \Vert^2\\
				&  + \frac{1-\tau}{2\rho}\big( \Vert \boldlambda^{k+1}-\boldlambda^k\Vert^2 + \Vert \boldlambda^{k+1}\Vert^2 - \Vert \boldlambda^k \Vert^2\big)  \\
				&   +\frac{ \rho}{2}\Vert \A\x^{k +1} - \bb \Vert^2 +  c \Big(  \frac{1-2\tau^2}{2\rho}  \Vert \boldlambda^{k+1}\! -\!\boldlambda^k \Vert^2  \\
				& + {1}/{2}\Vert  \x^{k+1}-\x^k\Vert^2_{\Q} +{L_g}/{2}\Vert \x^{k} - \x^{k-1} \Vert^{2} \Big)  \\
			\end{split}
		\end{equation}
	}

  By combing \eqref{eq:sufficient_decrease} and \eqref{eq:T_def2}, we have 
	{\small 
		\setlength{\abovedisplayskip}{-2pt}
		\setlength{\belowdisplayskip}{1pt}
		\begin{align}\label{eq:lambda-inequality}
			\frac{1-\tau}{2\rho} \big(  \Vert \boldlambda^{k+1}\Vert^2 - \Vert \boldlambda^k \Vert^2\big)  + \frac{\tau}{\rho}\Vert \boldlambda^{k + 1}\Vert^2 \leq T_c^0. 
		\end{align}
	}
The above holds because we have  $f(\x) \geq 0$, $g(\x) \geq 0$ over $\X$ and the other terms are all non-negative.

	We next prove $\frac{\tau}{ \rho} \Vert \boldlambda^{k + 1}\Vert^2 \leq T_c^{0}$ by induction. For $k =0$, we can properly pick the initial point to satisfy the inequality. For iteration $k$, we assume $\frac{\tau}{ \rho} \Vert \boldlambda^{k}\Vert^2 \leq T_c^{0}$. We consider the two possible cases for iteration $k+1$, i.e., if $\Vert \boldlambda^{k+1} \Vert^2 \leq \Vert \boldlambda^k \Vert^2$, we straightforwardly have $\frac{\tau}{\rho} \Vert \boldlambda^{k +1}\Vert^2 \leq \frac{\tau}{\rho} \Vert \boldlambda^{k}\Vert^2 \leq T_c^0$, and else if $\Vert \boldlambda^{k+1} \Vert^2 \geq \Vert \boldlambda^k \Vert^2$,  we also have $\frac{\tau}{\rho} \Vert \boldlambda^{k +1}\Vert^2 \leq T_c^0$ by \eqref{eq:lambda-inequality}. 
	We therefore conclude $\Vert \boldlambda^{*} \Vert^2 \leq  \rho \tau^{-1} T_c^0$.

ii) Prove $T_c^0 \leq  \big( 4 + c(1 - 2\tau^2) + c/2\big) d_F+ {c L_g}/{2}\Vert \x^0 \Vert^2 +  {c \rho_F}/{4}~d_{\x}$: Invoke \textbf{Prop.} \ref{prop:prop-L} and set $k = 0$, we have 
\begin{align*}
	& \Lag_{\rho}(\x^1, \boldlambda^1) - \frac{\tau}{2 \rho}\Vert \boldlambda^1 \Vert^2  \leq 	\Lag_{\rho}(\x^0, \boldlambda^0) - \frac{\tau}{2 \rho}\Vert \boldlambda^0 \Vert^2 \\
	&  - \Vert \x^{1} - \x^0 \Vert_{\Q}^2 + \frac{\rho_F}{2} \Vert \x^1-\x^{0} \Vert^2 -\frac{\rho}{2} \Vert \A(\x^{1} - \x^0)\Vert^2  \\
	& + \frac{2-\tau}{2\rho} \Vert \boldlambda^{1} - \boldlambda^0\Vert^2.
\end{align*}

By invoking \eqref{eq:second-term_L} and setting $\boldlambda^{-1} = 0 $,  $\boldlambda^0 =0$, $\A \x^0 = b$, $\Q: =  \rho G_{\A} + \beta G_{\B} - \rho \A^\top\A $, we have 
\begin{align*}
	& \frac{\rho}{2}\Vert \A \x^1 - \bb \Vert^2+ \frac{2\Q + \rho \A^\top\A -\rho_f \mathbf{I}_N }{2} \Vert \x^1 - \x^0 \Vert^2\\
	&  \leq f(\x^0) + g(\x^0) -  f(\x^1) - g(\x^1) \\
\end{align*}
Since we have $f(\x) \geq 0$ and $f(\x) \geq 0$ over the set $\X$, we have (the  term $\frac{\rho\A^\top\A}{2}\Vert \x^1 - \x^0 \Vert^2$ is non-negative)
\begin{align}
	& \label{eq:C4}\frac{\rho}{2}\Vert \A \x^1 - \bb \Vert^2 \leq d_F. \\
	&  \frac{2 \Q - \rho_F \mathbf{I}_N}{2}\Vert \x^1 - \x^0 \Vert^2 \leq d_F  \notag\\
	\label{eq:C5}&   \Vert \x^1 - \x^0\Vert^2_{\Q} \leq d_F + \rho_F/2 d_{\x}. 
\end{align}
where the last inequality is by $d_{\x}:= \max_{\x, \y} \Vert \x - \y \Vert^2$.

Further, based on the dual update, we have
\begin{align}
	 \label{eq:C6}& \frac{1}{2\rho}\Vert  \boldlambda^1\Vert^2 = \frac{\rho}{2}\Vert \A \x^1 - \bb\Vert^2 \leq d_F
\end{align}

Further, we have
{\small 
\begin{align*}
T_c^0&  = f(\x^{1}) + g(\x^{1})  + \frac{2 + c(1 - 2\tau^2)}{2\rho} \Vert \boldlambda^{1} \Vert^2 + \frac{\rho}{2}\Vert \A\x^{1} - \bb \Vert^2  \notag\\
	&  + \frac{c}{2}\Vert  \x^{1}-\x^0\Vert^2_{\Q} +\frac{cL_g}{2}\Vert \x^{0} \Vert^{2} ~~\text{by \eqref{eq:T_def2}~and $\boldlambda^0 = 0$}   \notag\\
	\leq &~~d_F + (2 + c(1 - 2\tau^2)) d_F + d_F \\
	& \quad \quad +\frac{c}{2} d_F + \frac{c \rho_F}{4}d_{\x}+ \frac{cL_g}{2}\Vert \x^0\Vert^2 ~~\text{by \eqref{eq:C4}, \eqref{eq:C5}, \eqref{eq:C6}}\\
	=  &~~ \big( 4 + c(1 - 2\tau^2) + c/2\big) d_F+ {c L_g}/{2}\Vert \x^0 \Vert^2 +  {c \rho_F}/{4}~d_{\x}
\end{align*}
}

By combining  i) and ii), we have  $\tau^2 \rho^{-2} \Vert \boldlambda^{*} \Vert^2 \leq \epsilon^2$. By invoking \textbf{Theorem} \ref{thm:theorem1}, we directly have
	\begin{align*}
	& \text{dist}\big( \nabla  f(\starx) + \nabla g(\starx) + \A^\top \hat{\boldlambda}^{*} + N_{\X}(\starx) , 0\big) \\
	&+ \Vert \A \starx - \bb  \Vert  \leq  \tau \rho^{-1} \Vert \boldlambda^{*} \Vert \leq \epsilon,  
\end{align*} 
which closes the proof. 
\end{document}